\documentclass[journal]{IEEEtran}
\usepackage{bm}
\usepackage{hyperref}
\usepackage{float}
\usepackage{multirow}
\usepackage{tikz}
\usetikzlibrary{arrows.meta,calc,fit,positioning,quotes}
\usepackage{booktabs}
\usepackage{threeparttable}

\usepackage{amssymb}
\usepackage{amsmath}
\usepackage{graphicx}
\usepackage[noend]{algpseudocode}
\usepackage{algorithm}
\usepackage{colortbl}
\usepackage[numbers,sort&compress]{natbib}
\usepackage{hyperref}
\newtheorem{lemma}{Lemma}[section]
\newtheorem{definition}{Definition}[section]
\newtheorem{theorem}{Theorem}[section]
\newtheorem{corollary}{Corollary}[section]
\newtheorem{example}{Example}[section]
\newtheorem{proposition}{Proposition}[section]
\newtheorem{property}{Property}[section]
\newtheorem{application}{Application}[section]

\begin{document}
%\title{Individual convergence of non-convex stochastic momentum methods in deep learning}
\title{Convex Quaternion Optimization for Signal Processing: Theory and Applications}

\author{Shuning Sun, Qiankun Diao, Dongpo Xu, Pauline Bourigault and Danilo P. Mandic,~\IEEEmembership{Fellow,~IEEE} 

%\author{Jinlan Liu, Yinghua Lu, Dongpo Xu, Jun Kong, Miao Qi and Danilo P. Mandic, \textit{Fellow}, \textit{IEEE}% <-this % stops a space

\thanks{This work was funded in part by the National Natural Science Foundation of China (No. 62176051), in part by National Key R$\&$D Program of China (No. 2021YFA1003400), and in part by the Fundamental Research Funds for the Central Universities of China (No. 2412020FZ024). (\textit{Corresponding author: Dongpo Xu})}% <-this % stops a space
\thanks{Shuning Sun, Qiankun Diao and Dongpo Xu are with the Key Laboratory for Applied Statistics of MOE, School of Mathematics and Statistics, Northeast Normal University, Changchun, 130024, China. (e-mail: xudp100@nenu.edu.cn)}
\thanks{Pauline Bourigault and Danilo P. Mandic are with the Department of Electrical and Electronic Engineering,
Imperial College London, London SW7 2AZ, UK. (e-mail: d.mandic@imperial.ac.uk)}}

%(\textit{Corresponding authors: Dongpo Xu, Yinghua Lu})
%(e-mail: dongpoxu@gmail.com)


\maketitle

% in the abstractThe average
\begin{abstract}
	Convex optimization methods have been extensively used in the fields of communications and signal processing. However, the theory of quaternion optimization is currently not as fully developed and systematic as that of complex and real optimization.
	To this end, we establish an essential theory of convex quaternion optimization for signal processing based on the generalized Hamilton-real (GHR) calculus. This is achieved in a way which conforms with traditional complex and real optimization theory. For rigorous, We present five discriminant theorems for convex quaternion functions, and four discriminant criteria for strongly convex quaternion functions. Furthermore, we provide a fundamental theorem for the optimality of convex quaternion optimization problems, and demonstrate its utility through three applications in quaternion signal processing. These results provide a solid theoretical foundation for convex quaternion optimization and open avenues for further developments in signal processing applications. 
\end{abstract}


% \fi
\IEEEpeerreviewmaketitle

\section{Introduction}
\IEEEPARstart{Q}{uaternions} were first introduced by William Hamilton in 1843 as an associative but non-commutative algebra over the real numbers \cite{hamilton1840new}. Since then, they have become a powerful tool in many fields, including 
%robotics \cite{10065552}, 
image processing \cite{qi2022quaternion,9782722}, signal processing \cite{took2011augmented,9001262,5447699}, and machine learning \cite{10048541,walia2023unveiling,2011Quaternion}.
Examples include the work by Jia et al. \cite{9782722}, who introduced a robust method for quaternion matrix completion, that can be used to reconstruct large-scale color images. 
Flamant et al. \cite{flamant2019time} demonstrated the efficiency of Quaternion Fourier Transform (QFT) in processing bivariate signals. 
Ogunfunmi et al. \cite{8853311} presented a kernel adaptive filter for quaternion data.
Moreover, Meng{\"u}{\c{c}} et al. \cite{mengucc2020design}
designed quaternion-valued second-order Volterra adaptive filters for nonlinear 3-D and 4-D signal processing. 
%Xiang et al. \cite{xiang2019performance} 
%analyzed the mean and mean square convergence of the deficient length quaternion least mean square  algorithms. 
Xia et al. \cite{8754780} established an estimation framework for processing quaternion-valued Gaussian data.
Finally, Zhang et al. \cite{10048541} discussed a new method for reducing the computation cost of quaternion signal estimation. 
Enshaeifar et al. \cite{008897c057f946ae922fb3af210efd39} introduced quaternion-valued singular spectrum analysis for multichannel electroencephalogram analysis.

The theory of real-valued and complex-valued convex optimization is well-established and has seen widely used in the areas of communications \cite{luo2006introduction}, machine learning \cite{sra2012optimization,jaggi2011sparse,krejic2023inexact} and signal processing \cite{7954619,8008828,gershman2010convex}. In recent years, convex quaternion optimization has also attracted interest. For example, 
Qi et al. \cite{qi2022quaternion} studied first-order derivatives and second-order partial derivatives of real-valued functions of quaternion variables over their real and imaginary $i$, $j$, $k$ parts. However, this complicates the proof and computational process in quaternion optimization. Flamant et al. \cite{flamant2021general} and Liu et al. \cite{liu2019constrained} provided first-order characterization of quaternion functions by generalized Hamilton-real (GHR) calculus \cite{xuRSOS}. However, these useful attempts lack the discussion of 
%the ??
gradient monotonicity and second-order characterization for convex quaternion function, a pre-requisite for practical applications. 

To fill this void, we have systematically address the theory of convex optimization in the quaternion domain. For rigorous, this is achieved based on the GHR calculus \cite{xuRSOS}, a generalization of Wirtinger-calculus \cite{hjorungnes2011complex,wirtinger1927formalen,brandwood1983complex} from the complex domain to the quaternion field. 
Before the introduction of the GHR calculus, the quaternion pseudo-derivative was used for calculating the gradient, which transforms the quaternion optimization problem into a lengthy and complicated real optimization problem;  the solution is then found by using real-valued optimization algorithms \cite{arena1997multilayer,yoshida2004model}. Flamant et al. \cite{9001262,flamant2021general} demostrated that the GHR calculus is a powerful theory in quaternion signal processing and non-negative matrix factorization. Meng{\"u}{\c{c}} et al. \cite{mengucc2018novel,mengucc2020design} established that the GHR calculus paves the way for the theory and applications of quaternion-valued adaptive filters. 
Took and Xia \cite{8682795} proposed a multichannel quaternion least-mean-square based on the GHR calculus for the adaptive filtering.
Parcollet et al. \cite{parcollet2020survey} further emphasized the significance of the GHR calculus as a recent breakthrough in the field.

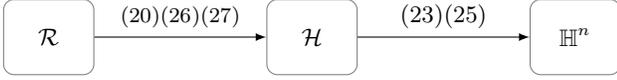
\begin{figure}[H]
	\begin{center}
		%\includegraphics[width=0.7\linewidth]{WechatIMG352.png}
		\begin{tikzpicture}[
			squarednode1/.style={rectangle, rounded corners, draw=black!50, fill=white!5, minimum width=1.2cm, minimum height=1.0cm},
			squarednode2/.style={rectangle, rounded corners, draw=blue!50, fill=green!5, minimum width=1.2cm, minimum height=0.8cm},font=\small,align=center]
			%\draw[help lines, color=gray] (-5,-5) grid (5,5);
			\node[squarednode1] (h1) at (0,0) { 
				$\mathcal{H}$} ;
			%\node[squarednode2] (h2) at (2.8,1.2) {\ref{BG}} ; 
			\node[squarednode1] (h3) at (3.5,0) { 
				$\mathbb{H}^n$} ;
			%\node[squarednode2] (h3) at (2.8,0) {88\ref{BV}} ; 
			%\node[squarednode2] (h4) at (0,-2.1) {\ref{ES}} ; 
			\node[squarednode1] (h5) at (-3.5,0) { 
				$\mathcal{R}$} ;
			%\node[squarednode2] (h5) at (-2.8,0) {\ref{ABC}} ; 
			%\node[squarednode2] (h6) at (-2.8,-1.2) {\ref{SG}} ; 
			%\node[squarednode2] (h7) at (0,2.1) {\ref{SS}} ; 
			%\draw[-latex] (h2) -- node[above,sloped] {}(h1) ; 
			\draw[-latex] (h1) -- node[above,sloped]  {$(\ref{h})(\ref{HHH})$}(h3) ; 
			%\draw[-latex] (h4) -- node[above,right]  {\footnotesize PL}(h1) ; 
			%\draw[-latex] (h6) -- node[above,sloped] {}(h1);
			\draw[-latex] (h5) -- node[above,sloped] {\footnotesize $(\ref{r,h})(\ref{grad})(\ref{hassion})$}(h1);  
			%\draw[-latex] (h7) -- node[above,right]{\footnotesize PL}(h1); 
		\end{tikzpicture}
		% The figure caption is below the figure.
	\end{center}
	%\vspace{-1 in}
	\caption{The derivation process of quaternion optimization theory, with the sets $\mathcal{R}$ and $\mathcal{H}$ defined by $(\ref{r4N}),(\ref{calH})$.}
	\label{pppp}
\end{figure}

The theory of convex optimization in the quaternion field has gained attention due to its promising applications in signal processing and optimization. The aim of this work is to  develop the convexity theory of quaternion function using the GHR calculus \cite{xuRSOS}. To this end, we make use of the duality of the augmented quaternion vector 
$\bm{q}_{\mathcal{H}} \triangleq \left(\bm{q}^{\mathsf{T}},\bm{q}^{i\mathsf{T}},\bm{q}^{j\mathsf{T}},\bm{q}^{k\mathsf{T}}\right)^{\mathsf{T}}$ and the augmented real vector $\bm{q}_{\mathcal{R}} \triangleq \left(\bm{q}_a^{\mathsf{T}},\bm{q}_b^{\mathsf{T}},\bm{q}_c^{\mathsf{T}},\bm{q}_d^{\mathsf{T}}\right)^{\mathsf{T}}$ 
\cite{xuTNNLS}.
%, which are related by (\ref{r,h}). 
Next, we employ the relationships between augmented quaternion gradient and augmented real gradient, as well as between the augmented quaternion Hessian matrix and augmented real Hessian matrix, as shown in \cite{xuTNNLS}.
%, and in (\ref{grad})  (\ref{hassion}). 
Based on these results, we extend the discriminant criteria for convexity from the real field to the augmented quaternion space, $\mathcal{H}$, and then to the quaternion field $\mathbb{H}^{n}$, as illustrated in Figure $\ref{pppp}$. Moreover, we define and present four discriminant theorems for strong convexity, by employing the discriminant criteria of convex quaternion functions. Finally, we present a fundamental theorem for the optimality of convex quaternion problems and provide three illustrative applications in the field of signal processing, including quaternion linear mean-square error filter, quaternion projection on affine equality constraint, and quaternion minimum variance beamforming.
This work makes three significant contributions to the theory of convex optimization in the quaternion field: 
%\vspace{-\topsep}
\begin{itemize}
	\item By using the GHR calculus, we establish five discriminant theorems for convex functions in the quaternion field. These theorems include gradient monotonicity and second-order characterization. 
	\item We provide a clear definition and four discriminant criteria for strongly convex functions in the quaternion field; these are consistent with their counterpart real and complex convexity theorems. 
	%These results provide a comprehensive understanding of strongly convex function behavior. 
	\item A fundamental theorem is proposed for the optimality of convex quaternion problems, together with some practical applications of convex quaternion optimization in communications and signal processing.
\end{itemize}

This paper is organized as follows. In Section II, we give an overview of quaternion algebra, the GHR calculus, and some equivalence relationships. Section III presents five discriminant theorems for convex quaternion functions, covering first-order characterization, second-order characterization and some examples of convex quaternion functions. Section IV introduces the definition and discriminant theorems for strongly convex quaternion functions. In Section V, we propose a fundamental theorem for convex quaternion optimization problems and provide three practical applications in signal processing. Finally, this paper concludes with Section VI.

\section{Preliminaries}
\subsection{Quaternion Algebra}
A quaternion, $q$, can be expressed as
\begin{equation}
	q=q_a+q_bi+q_cj+q_dk,
\end{equation}
where $q_a, q_b,q_c,q_d\in \mathbb{R}$,  and the imaginary
units $i$, $j$ and $k$ satisfy $i^2=j^2=k^2=ijk=-1$, $ij=-ji=k$, $jk=-kj=i$, $ki=-ik=j$. 
The set of quaternions is defined as 
$\mathbb{H} \triangleq \left\{q=q_a+q_bi+q_cj+q_dk \mid q_a,q_b,q_c,q_d\in \mathbb{R}\right\}$. Owing to the properties of the imaginary units, 
the multiplication of
two quaternions in $\mathbb{H}$ is noncommutative. 
The real part of $q$ is denoted by ${\rm Re}\left\{q\right\}=q_a$, whereas the imaginary part (pure quaternion) is ${\rm Im}\left\{q\right\}=q_bi+q_cj+q_dk$.
The conjugate of q is $q^* = {\rm Re}\left\{q\right\} - {\rm Im}\left\{q\right\} = q_a-q_bi-q_cj-q_dk$. %, and the conjugate of the product satisfies $(pq)^*=q^*p^*$, $\forall p,q \in \mathbb{H}$.
The modulus of a quaternion is defined as $\mid q\mid =\sqrt{qq^*}$. Define also $q^\mu$, as it is used in Definition $\ref{2.1}$.
\begin{definition}
	[Quaternion rotation \cite{ward2012quaternions}]\label{2.1}
	For any quaternion, $q$, and a nonzero quaternion $\mu$, the transformation
	\begin{equation}\label{xz}
		q^\mu \triangleq \mu q \mu^{-1}
	\end{equation}
	describes a rotation of $q$.
%	 by an angle 2$\theta$ about the vector part of $\mu$, where $\mu = \mid\mu\mid({\rm cos}\theta+\hat{\mu}{\rm sin}\theta)$
%	is any non-zero quaternion, $\hat{\mu}=\frac{V_\mu}{\mid V_\mu \mid}$, $\theta={\rm arccos}\left(\frac{S_\mu}{\mid \mu\mid }\right)$.
	
	In particular, if $\mu$ in $(\ref{xz})$ is a pure unit quaternion, then the
	quaternion rotation in $(\ref{xz})$ becomes quaternion involution \cite{ell2007quaternion}, such as
	\begin{align}
		q^{i} & =-i q i=q_{a}+i q_{b}-j q_{c}-k q_{d}, \\
		q^{j} & =-j q j=q_{a}-i q_{b}+j q_{c}-k q_{d}, \\
		q^{k} & =-k q k=q_{a}-i q_{b}-j q_{c}+k q_{d}.
	\end{align}
\end{definition}
\begin{property}
	[Properties of quaternion rotation \cite{Mandic2015The}] For any 
	$p, q \in \mathbb{H}$, and $\forall \nu, \mu \in \mathbb{H}$, the following holds
	\begin{equation}\label{xzp}
		\begin{array}{c}
			(p q)^{\mu}=p^{\mu} q^{\mu}, ~
			q^{\mu \nu}=\left(q^{\nu}\right)^{\mu},\\
			q^{\mu *} \triangleq \left(q^{*}\right)^{\mu} = \left(q^{\mu}\right)^{*} \triangleq q^{* \mu}.
		\end{array}
	\end{equation}
\end{property}
\subsection{The GHR Calculus}
\begin{definition}
	[Real-differentiability \cite{skw}]
	A quaternion function $f:\mathbb{H} \to \mathbb{H}$, given by  
	$ f(q)=f_{a}\left(q_{a}, q_{b}, q_{c}, q_{d}\right)+i f_{b}\left(q_{a}, q_{b}, q_{c}, q_{d}\right)+ j f_{c}\left(q_{a}, q_{b}, q_{c}, q_{d}\right)+k f_{d}\left(q_{a}, q_{b}, q_{c}, q_{d}\right) $ is called real differentiable, if $f_{a}$, $f_{b}$, $f_{c}$, $f_{d}$ are differentiable as functions of the real variables $q_{a}$, $q_{b}$, $q_{c}$, $q_{d}$.
\end{definition}
\begin{definition}
	[GHR derivatives \cite{xuRSOS}] \label{GHR}
	If $f: \mathbb{H} \rightarrow \mathbb{H}$ is real-differentiable, then the left GHR
	derivatives of the function $f$ with respect to $q^{\mu}$ and $q^{\mu *}$ $\left(\mu \neq 0,\mu \in \mathbb{H}\right)$ are defined as
	\begin{align}
		\dfrac{\partial f}{\partial q^{\mu}} &=\dfrac{1}{4}\left(\dfrac{\partial f}{\partial q_{a}}-\dfrac{\partial f}{\partial q_{b}} i^{\mu}-\dfrac{\partial f}{\partial q_{c}} j^{\mu}-\dfrac{\partial f}{\partial q_{d}} k^{\mu}\right) \in \mathbb{H},\\
		\dfrac{\partial f}{\partial q^{\mu *}} &=\dfrac{1}{4}\left(\dfrac{\partial f}{\partial q_{a}}+\dfrac{\partial f}{\partial q_{b}} i^{\mu}+\dfrac{\partial f}{\partial q_{c}} j^{\mu}+\dfrac{\partial f}{\partial q_{d}} k^{\mu}\right) \in \mathbb{H},
	\end{align}
	where $q=q_a+q_bi+q_cj+q_dk$, $q_a, q_b,q_c,q_d\in \mathbb{R}$, and $\frac{\partial f}{\partial q_{a}},\frac{\partial f}{\partial q_{b}},\frac{\partial f}{\partial q_{c}},\frac{\partial f}{\partial q_{d}} \in \mathbb{R}$ are the partial derivatives
	of $f$ with respect to $q_{a}$, $q_{b}$, $q_{c}$, $q_{d}$.
\end{definition}
\begin{property}
	[Properties of the GHR derivatives \cite{xuRSOS}]
	If $f: \mathbb{H} \rightarrow \mathbb{H}$, $g: \mathbb{H} \rightarrow \mathbb{H}$, then\\
	Product rule:
	\begin{equation}\label{pr}
		\dfrac{\partial(f g)}{\partial q^{\mu}}=f \dfrac{\partial g}{\partial q^{\mu}}+\dfrac{\partial f}{\partial q^{g \mu}} g, \quad
		\dfrac{\partial(f g)}{\partial q^{\mu *}}=f \dfrac{\partial g}{\partial q^{\mu *}}+\dfrac{\partial f}{\partial q^{g \mu *}} g
	\end{equation}
	\\Chain rule:
	\begin{align}
		\dfrac{\partial f(g(q))}{\partial q^{\mu}}&=\sum_{\nu \in\{1, i, j, k\}} \dfrac{\partial f}{\partial g^{\nu}} \dfrac{\partial g^{\nu}}{\partial q^{\mu}},\\
		\dfrac{\partial f(g(q))}{\partial q^{\mu *}}&=\sum_{\nu \in\{1, i, j, k\}} \dfrac{\partial f}{\partial g^{\nu *}} \dfrac{\partial g^{\nu *}}{\partial q^{\mu *}}
	\end{align}
	\\Rotation rule:
	\begin{equation}\label{rotation}
		\left(\dfrac{\partial f}{\partial q^{\mu}}\right)^{\nu}
		=
		\dfrac{\partial f^{\nu}}{\partial q^{\nu \mu}}, 
		\quad\left(\dfrac{\partial f}{\partial q^{\mu *}}\right)^{\nu}
		=
		\dfrac{\partial f^{\nu}}{\partial q^{\nu \mu *}}
	\end{equation}
	Conjugate rule: 
	%	\begin{equation}
		%		\left(\dfrac{\partial f}{\partial q^{\mu}}\right)^{*}
		%		=
		%		\dfrac{\partial f^{*}}{\partial q^{\mu *}},
		%		\quad\left(\dfrac{\partial f}{\partial q^{\mu *}}\right)^{*}
		%		=
		%		\dfrac{\partial f^{*}}{\partial q^{\mu}}
		%	\end{equation}
	If $f: \mathbb{H} \rightarrow \mathbb{R}$, 
	\begin{equation}\label{cong}
		\left(\dfrac{\partial f}{\partial q^{\mu}}\right)^{*}
		=
		\dfrac{\partial f}{\partial q^{\mu *}},
		\quad\left(\dfrac{\partial f}{\partial q^{\mu *}}\right)^{*}
		=
		\dfrac{\partial f}{\partial q^{\mu}}.
	\end{equation}
\end{property}
\begin{definition}
	[Quaternion gradient \cite{xuTNNLS}]
	The quaternion
	gradient and its conjugate gradient of a function $f: \mathbb{H}^n \rightarrow \mathbb{R} $ are defined as
	\begin{equation}
		\nabla_{\bm{q}} f 
		\triangleq
		\left(\dfrac{\partial f}{\partial \bm{q}}\right)^{\mathsf{T}}
		=
		\left(\dfrac{\partial f}{\partial q_{1}}, \ldots, \dfrac{\partial f}{\partial q_{n}}\right)^{\mathsf{T}} \in \mathbb{H}^n,
	\end{equation}
	\begin{equation}\label{getd}
		\nabla_{\bm{q}^{*}} f \triangleq\left(\dfrac{\partial f}{\partial \bm{q}^{*}}\right)^{\mathsf{T}}=\left(\dfrac{\partial f}{\partial q_{1}^{*}}, \ldots, \dfrac{\partial f}{\partial q_{n}^{*}}\right)^{\mathsf{T}} \in \mathbb{H}^n,
	\end{equation}
	where $\left(\frac{\partial f}{\partial \bm{q}}\right)^{\mathsf{T}}$ is the transpose of $\frac{\partial f}{\partial \bm{q}}$.
\end{definition}
\begin{definition}
	[Quaternion Hessian \cite{xuTNNLS}]
	Let $f: \mathbb{H}^n \rightarrow \mathbb{R} $, then the two quaternion Hessian matrices are defined as
	\begin{align}
		\boldsymbol{H}_{\bm{q} \bm{q}} &\triangleq \dfrac{\partial}{\partial \bm{q}}\left(\dfrac{\partial f}{\partial \bm{q}}\right)^{\mathsf{T}}\nonumber\\
		&=\left(\begin{array}{ccc}
			\dfrac{\partial^{2} f}{\partial q_{1} \partial q_{1}} & \cdots & \dfrac{\partial^{2} f}{\partial q_{n} \partial q_{1}} \\
			\vdots & \ddots & \vdots \\
			\dfrac{\partial^{2} f}{\partial q_{1} \partial q_{n}} & \cdots & \dfrac{\partial^{2} f}{\partial q_{n} \partial q_{n}}
		\end{array}\right) \in \mathbb{H}^{n\times n},
	\end{align}
	\begin{align}	\label{hqq*}
		\boldsymbol{H}_{\bm{q} \bm{q}^{*}} 
		&\triangleq 
		\dfrac{\partial}{\partial \bm{q}}\left(\dfrac{\partial f}{\partial \bm{q}^{*}}\right)^{\mathsf{T}}\nonumber\\
		&=\left(\begin{array}{ccc}
			\dfrac{\partial^{2} f}{\partial q_{1} \partial q_{1}^{*}} & \cdots & \dfrac{\partial^{2} f}{\partial q_{n} \partial q_{1}^{*}} \\
			\vdots & \ddots & \vdots \\
			\dfrac{\partial^{2} f}{\partial q_{1} \partial q_{n}^{*}} & \cdots & \dfrac{\partial^{2} f}{\partial q_{n} \partial q_{n}^{*}}
		\end{array}\right) \in \mathbb{H}^{n\times n}.
	\end{align}
\end{definition}
\subsection{The Relationship of Augmented Quaternion and the Augmented Real Vector, Gradient, and Hessian Matrix}
Consider a quaternion vector $\bm{q} = \bm{q}_a + \bm{q}_b {i} + \bm{q}_c {j} + \bm{q}_d {k} \in \mathbb{H}^n$ where $\bm{q}_a,\bm{q}_b,\bm{q}_c,\bm{q}_d\in\mathbb{R}^n$. Define its augmented real vector as $\bm{q}_{\mathcal{R}} \triangleq \left(\bm{q}_a^{\mathsf{T}},\bm{q}_b^{\mathsf{T}},\bm{q}_c^{\mathsf{T}},\bm{q}_d^{\mathsf{T}}\right)^{\mathsf{T}} \in \mathcal{R}$ \cite{took2011augmented,2010Properness} and the augmented quaternion vector as
$\bm{q}_{\mathcal{H}} \triangleq \left(\bm{q}^{\mathsf{T}},\bm{q}^{i\mathsf{T}},\bm{q}^{j\mathsf{T}},\bm{q}^{k\mathsf{T}}\right)^{\mathsf{T}} \in \mathcal{H}$ \cite{xuTNNLS}, where the set of augmented real vectors and the set of augmented quaternion vectors are defined as
\begin{align}
	\mathcal{R} &\triangleq \left\{\bm{q}_{\mathcal{R}}=\left(\bm{q}_a^{\mathsf{T}},\bm{q}_b^{\mathsf{T}},\bm{q}_c^{\mathsf{T}},\bm{q}_d^{\mathsf{T}}\right)^{\mathsf{T}} \mid  \bm{q} \in \mathbb{H}^n\right\} = \mathbb{R}^{4n},\label{r4N}\\
	\mathcal{H} &\triangleq \left\{\bm{q}_{\mathcal{H}}= \left(\bm{q}^{\mathsf{T}},\bm{q}^{i\mathsf{T}},\bm{q}^{j\mathsf{T}},\bm{q}^{k\mathsf{T}}\right)^{\mathsf{T}} \mid \bm{q}\in\mathbb{H}^n\right\} \subset \mathbb{H}^{4n}.\label{calH}
\end{align}
By definition, there exists a one-to-one mapping between $\mathbb{H}^n$, $\mathcal{R}$ and $\mathcal{H}$ \cite{flamant2021general}.
\begin{proposition}[\cite{xuTNNLS}]
	The relationship between the augmented quaternion vector, $\bm{q}_{\mathcal{H}}$, and the augmented real vector, $\bm{q}_{\mathcal{R}}$, is given by
	\begin{equation}\label{r,h}
		\bm{q}_{\mathcal{H}} = \bm{J}_n \bm{q}_{\mathcal{R}} 
		\quad \Leftrightarrow \quad
		\bm{q}_{\mathcal{R}} = \dfrac{1}{4} \bm{J}_n^{\mathsf{H}} \bm{q}_{\mathcal{H}},
	\end{equation}
	where
	\begin{align}
		\bm{J}_n = \left(\begin{array}{cccc}
			\bm{I}_n & i \bm{I}_n & j \bm{I}_n & k \bm{I}_n \\
			\bm{I}_n & i \bm{I}_n & -j \bm{I}_n & -k \bm{I}_n \\
			\bm{I}_n & -i \bm{I}_n & j \bm{I}_n & -k \bm{I}_n \\
			\bm{I}_n & -i \bm{I}_n & -j \bm{I}_n & k \bm{I}_n
		\end{array}\right)
		\in \mathbb{H}^{4n \times 4n},
	\end{align}
	and $\bm{J}_n^{\mathsf{H}} \bm{J}_n = 4\bm{I}_{4n}$, while $\bm{I}_n$ is the $n \times n$ identity matrix, with $\bm{J}_n^{\mathsf{H}}$ as the conjugate transpose of $\bm{J}_n$.  
\end{proposition}

From $(\ref{r,h})$, the quaternion function $f\left(\bm{q}\right):\mathbb{H}^n\rightarrow\mathbb{R}$ can be
viewed in three equivalent forms
\cite{xuTNNLS}, as follows
\begin{align} \label{iff}
	\begin{aligned}
		f\left(\bm{q}\right) &\quad \Leftrightarrow \quad f\left(\bm{q}_{\mathcal{R}}\right)
		\triangleq f\left(\bm{q}_a,\bm{q}_b,\bm{q}_c,\bm{q}_d\right)\\ &\quad \Leftrightarrow \quad f\left(\bm{q}_{\mathcal{H}}\right) \triangleq  f\left(\bm{q},\bm{q}^i,\bm{q}^j,\bm{q}^k\right).
	\end{aligned}
\end{align}
Note that these three functions are equivalent but have different forms, denoted as $f$ for simplicity. 
Here, the variables of the functions $f\left(\bm{q}\right)$, $f\left(\bm{q}_{\mathcal{R}}\right)$, and $f\left(\bm{q}_{\mathcal{H}}\right)$ are quaternion vectors, augmented real vectors, and augmented quaternion vectors, respectively. They are referred to as quaternion function, augmented real function, and augmented quaternion function, respectively. 
For $f\left(\bm{q}_{\mathcal{R}}\right):\mathcal{R}\rightarrow\mathbb{R}$, its augmented real gradient is defined as 
$\nabla_{\mathcal{R}} f \triangleq 
\left(\frac{\partial f}{\partial \bm{q}_{\mathcal{R}}}\right)^{\mathsf{T}}$ and the augmented real Hessian matrix as $\boldsymbol{H}_{\mathcal{R}\mathcal{R}} \triangleq \frac{\partial}{\partial \bm{q}_{\mathcal{R}}}\left(\frac{\partial f}{\partial \bm{q}_{\mathcal{R}}}\right)^{\mathsf{T}}$.
For   $f\left(\bm{q}_{\mathcal{H}}\right):\mathcal{H}\rightarrow\mathbb{R}$, the
augmented quaternion gradient and its conjugate gradient are defined as \cite{5623300}
\begin{align}\label{h}
	\nabla_{\mathcal{H}} f 
	&\triangleq 
	\left(\dfrac{\partial f}{\partial \bm{q}_{\mathcal{H}}}\right)^{\mathsf{T}}
	=
	\begin{pmatrix}
		\nabla _ {\bm{q}} f	\\
		\nabla _ {\bm{q}^i} f	\\
		\nabla _ {\bm{q}^j} f	\\
		\nabla _ {\bm{q}^k} f	
	\end{pmatrix},\\	
	\nabla_{\mathcal{H}^*} f 
	&\triangleq 
	\left(\dfrac{\partial f}{\partial \bm{q}_{\mathcal{H}}^*}\right)^{\mathsf{T}}
	=
	\begin{pmatrix}
		\nabla _ {\bm{q}^*} f	\\
		\nabla _ {\bm{q}^{i*}} f	\\
		\nabla _ {\bm{q}^{j*}} f	\\
		\nabla _ {\bm{q}^{k*}} f	
	\end{pmatrix},
\end{align}
and the augmented quaternion Hessian matrix is defined as 
\begin{align}\label{HHH}
	\boldsymbol{H}_{\mathcal{H}\mathcal{H}^*} &\triangleq 
	\dfrac{\partial}{\partial \bm{q}_{\mathcal{H}}}
	\left(\dfrac{\partial f}{\partial \bm{q}_{\mathcal{H}}^*}\right)^{\mathsf{T}}\\
	&=
	\left(\begin{array}{llll}
		\boldsymbol{H}_{\bm{q} \bm{q}^{*}} & \boldsymbol{H}_{\bm{q}^{i} \bm{q}^{*}} & \boldsymbol{H}_{\bm{q}^{j} \bm{q}^{*}} & \boldsymbol{H}_{\bm{q}^{k} \bm{q}^{*}} \\
		\boldsymbol{H}_{\bm{q} \bm{q}^{i *}} & \boldsymbol{H}_{\bm{q}^{i} \bm{q}^{i *}} & \boldsymbol{H}_{\bm{q}^{j} \bm{q}^{i *}} & \boldsymbol{H}_{\bm{q}^{k} \bm{q}^{i *}} \\
		\boldsymbol{H}_{\bm{q} \bm{q}^{j *}} & \boldsymbol{H}_{\bm{q}^{i} \bm{q}^{j *}} & \boldsymbol{H}_{\bm{q}^{j} \bm{q}^{j *}} & \boldsymbol{H}_{\bm{q}^{k} \bm{q}^{j *}} \\
		\boldsymbol{H}_{\bm{q} \bm{q}^{k *}} & \boldsymbol{H}_{\bm{q}^{i} \bm{q}^{k *}} & \boldsymbol{H}_{\bm{q}^{j} \bm{q}^{k *}} & \boldsymbol{H}_{\bm{q}^{k} \bm{q}^{k *}}
	\end{array}\right).\nonumber
\end{align}

\begin{proposition}
	[\cite{5623300}]
	The relationship between the augmented quaternion gradient, $\nabla_{\mathcal{H}^*} f$, and the augmented real gradient, $\nabla_{\mathcal{R}} f$, is given by 
	\begin{equation}\label{grad}
		\nabla_{\mathcal{H}^*} f = \dfrac{1}{4} \bm{J}_n \nabla_{\mathcal{R}} f
		\quad \Leftrightarrow \quad
		\nabla_{\mathcal{R}} f = {\bm{J}_n}^{\mathsf{H}} \nabla_{\mathcal{H}^*} f.
	\end{equation}
\end{proposition}

\begin{proposition}[\cite{xuTNNLS}]
	The relationship between the augmented quaternion Hessian matrix, $\boldsymbol{H}_{\mathcal{H}\mathcal{H}^*}$, and the augmented real Hessian matrix, $\boldsymbol{H}_{\mathcal{R}\mathcal{R}}$, is given by 
	\begin{equation}\label{hassion}
		\boldsymbol{H}_{\mathcal{H}\mathcal{H}^*} = \dfrac{1}{16} \bm{J}_n \boldsymbol{H}_{\mathcal{R}\mathcal{R}} \bm{J}_n^{\mathsf{H}}
		~ \Leftrightarrow ~
		\boldsymbol{H}_{\mathcal{R}\mathcal{R}} = \bm{J}_n^{\mathsf{H}} \boldsymbol{H}_{\mathcal{H}\mathcal{H}^*} \bm{J}_n
	\end{equation}
	where $\boldsymbol{H}_{\mathcal{H}\mathcal{H}^*} \triangleq 
	\frac{\partial}{\partial \bm{q}_{\mathcal{H}}}
	\left(\frac{\partial f}{\partial \bm{q}_{\mathcal{H}}^*}\right)^{\mathsf{T}}$, $\boldsymbol{H}_{\mathcal{R}\mathcal{R}} \triangleq \frac{\partial}{\partial \bm{q}_{\mathcal{R}}}\left(\frac{\partial f}{\partial \bm{q}_{\mathcal{R}}}\right)^{\mathsf{T}}$.
\end{proposition}

\begin{corollary}\label{hermite}
	The augmented quaternion Hessian matrix, 
	$\boldsymbol{H}_{\mathcal{H}\mathcal{H}^*}$, is a Hermite matrix, that is
	\begin{equation}\label{Hermite}
		\boldsymbol{H}_{\mathcal{H}\mathcal{H}^*}^{\mathsf{H}}
		=
		\boldsymbol{H}_{\mathcal{H}\mathcal{H}^*},
	\end{equation}
	where $\boldsymbol{H}_{\mathcal{H}\mathcal{H}^*}^{\mathsf{H}}$ is the conjugate transpose of $\boldsymbol{H}_{\mathcal{H}\mathcal{H}^*}$.
\end{corollary}
\begin{IEEEproof}
This is straightforward to demonstrate by using $(\ref{hassion})$ and the fact that $\boldsymbol{H}_{\mathcal{R}\mathcal{R}}$ is a Hermitian matrix.
%Since $\boldsymbol{H}_{\mathcal{R}\mathcal{R}}$ is a Hermite matrix and $(\ref{hassion})$, then 
%\begin{align}
%	\begin{aligned}
%		\boldsymbol{H}_{\mathcal{H}\mathcal{H}^*}^{\mathsf{H}} &= \dfrac{1}{16} \left(\bm{J}_n \boldsymbol{H}_{\mathcal{R}\mathcal{R}} \bm{J}_n^{\mathsf{H}}\right)^{\mathsf{H}}\\
%		&=
%		\dfrac{1}{16} \bm{J}_n \boldsymbol{H}_{\mathcal{R}\mathcal{R}}^{\mathsf{H}} \bm{J}_n^{\mathsf{H}}\\
%		&=
%		\dfrac{1}{16} \bm{J}_n \boldsymbol{H}_{\mathcal{R}\mathcal{R}} \bm{J}_n^{\mathsf{H}}\\
%		&=
%		\boldsymbol{H}_{\mathcal{H}\mathcal{H}^*}.
%	\end{aligned}
%\end{align}
%This completes the proof.
\end{IEEEproof}
\begin{proposition}
	For any 
	$\bm{p},\bm{q} \in \mathbb{H}^n$, their augmented real vectors are  $\bm{p}_{\mathcal{R}},\bm{q}_{\mathcal{R}}\in\mathcal{R}$, and their augmented quaternion vectors are $\bm{p}_{\mathcal{H}},\bm{q}_{\mathcal{H}}\in\mathcal{H}$. Then
	\begin{align}
		(a)&\quad \bm{p}_{\mathcal{H}}^{\mathsf{T}} \bm{q}_{\mathcal{H}}
		=
		4{\rm Re}\left\{\bm{p}^{\mathsf{T}} \bm{q}\right\};\label{pq}\\
		(b)&\quad 
		4\bm{p}_{\mathcal{R}}^{\mathsf{T}} \bm{q}_{\mathcal{R}}
		=
		\bm{p}_{\mathcal{H}}^{\mathsf{H}} \bm{q}_{\mathcal{H}}
		=
		4{\rm Re}\left\{\bm{p}^{\mathsf{H}} \bm{q}\right\} ;\label{pq2}\\
		(c)&\quad 2\lVert{\bm{p}_{\mathcal{R}}}\rVert_2
		=
		\lVert{\bm{p}_{\mathcal{H}}}\rVert_2
		=
		2\lVert{\bm{p}}\rVert_2; \label{fanshu}\\
		(d)&\quad \lVert{\bm{p}+\bm{q}}\rVert_2^2 
		=
		\lVert{\bm{p}}\rVert_2^2 + 2 {\rm Re}\left\{\bm{p}^{\mathsf{H}}\bm{q}\right\} +  \lVert{\bm{q}}\rVert_2^2.\label{sq}\quad
	\end{align}
\end{proposition}
\begin{IEEEproof} 	By the relationship of  $\bm{q}$, $\bm{q}_{\mathcal{R}}$, and $\bm{q}_{\mathcal{H}}$, we have 
\begin{align}
	&\begin{aligned}
		(a)~
		\bm{p}_{\mathcal{H}}^{\mathsf{T}} \bm{q}_{\mathcal{H}}
		&=
		\sum\limits_{\mu \in \left\{1,i,j,k \right\}} \bm{p}^{\mu \mathsf{T}} \bm{q}^{\mu} 
		\overset{\eqref{xzp}}{=}
		\sum\limits_{\mu \in \left\{1,i,j,k \right\}} \left(\bm{p}^\mathsf{T} \bm{q}\right)^{\mu}
		\\&\overset{\eqref{r,h}}{=}
		4{\rm Re}\left\{\bm{p}^{\mathsf{T}} \bm{q}\right\};
	\end{aligned}\\
	&\begin{aligned}
		(b)~
		4\bm{p}_{\mathcal{R}}^{\mathsf{T}} \bm{q}_{\mathcal{R}}
		&=
		4\bm{p}_{\mathcal{R}}^{\mathsf{H}} \bm{q}_{\mathcal{R}}
		\overset{\eqref{r,h}}{=}
		\bm{p}_{\mathcal{H}}^{\mathsf{H}} \bm{J}_n \dfrac{1}{4} \bm{J}_n^{\mathsf{H}} \bm{q}_{\mathcal{H}}
		=
		\bm{p}_{\mathcal{H}}^{\mathsf{H}} \bm{q}_{\mathcal{H}}
		\\&\overset{\eqref{pq}}{=}
		4{\rm Re}\left\{\bm{p}^{\mathsf{H}} \bm{q}\right\};
	\end{aligned}	\\
	&
	\begin{aligned}
		(c)~
		4\lVert{\bm{p}_{\mathcal{R}}}\rVert_2^2
		&=
		4\bm{p}_{\mathcal{R}}^{\mathsf{T}} \bm{p}_{\mathcal{R}}
		\overset{\eqref{pq2}}{=}
		\bm{p}_{\mathcal{H}}^{\mathsf{H}} \bm{p}_{\mathcal{H}}
		=
		\lVert{\bm{p}_{\mathcal{H}}}\rVert_2^2\\
		&\overset{\eqref{pq2}}{=}
		4{\rm Re}\left\{\bm{p}^{\mathsf{H}} \bm{p}\right\}
		=
		4\bm{p}^{\mathsf{H}} \bm{p}
		=
		4\lVert{\bm{p}}\rVert_2^2;
	\end{aligned}\\
	&
	\begin{aligned}
		(d)~ \lVert{\bm{p}+\bm{q}}\rVert_2^2 
		=
		&\left(\bm{p}+\bm{q}\right)^{\mathsf{H}}\left(\bm{p}+\bm{q}\right)
		\\=&
		\bm{p}^{\mathsf{H}}\bm{p} + \bm{p}^{\mathsf{H}}\bm{q} + \bm{q}^{\mathsf{H}}\bm{p} + \bm{q}^{\mathsf{H}}\bm{q}
		\\=&
		\lVert{\bm{p}}\rVert_2^2 + 2 {\rm Re}\left\{\bm{p}^{\mathsf{H}}\bm{q}\right\} +  \lVert{\bm{q}}\rVert_2^2.
	\end{aligned}
\end{align}
This completes the proof.
\end{IEEEproof}

\begin{proposition}
	If the quaternion function $f\left(\bm{q}\right):\mathbb{H}^n\rightarrow\mathbb{R}$ is real-differentiable, then $\forall \bm{p},\bm{q} \in \mathbb{H}^n$ we have
	\begin{align}
		&\begin{aligned}
			(a)~ \nabla_{\mathcal{R}} f\left(\bm{p}_{\mathcal{R}}\right)^{\mathsf{T}} \bm{q}_{\mathcal{R}}
			=&
			\nabla_{\mathcal{H}^*} f\left(\bm{p}_{\mathcal{H}}\right)^{\mathsf{H}} \bm{q}_{\mathcal{H}}
			\\=&
			4{\rm Re}\left\{\nabla_{\bm{p}^*} f\left(\bm{p}\right)^{\mathsf{H}} \bm{q}\right\};\label{tiduq}
		\end{aligned}\\
		&\begin{aligned}(b)~ 
			\nabla_{\mathcal{R}} f\left(\bm{p}_{\mathcal{R}}\right)^{\mathsf{T}} \nabla_{\mathcal{R}} &f\left(\bm{q}_{\mathcal{R}}\right)
			=
			4\nabla_{\mathcal{H}^*} f\left(\bm{p}_{\mathcal{H}}\right)^{\mathsf{H}} \nabla_{\mathcal{H}^*} f\left(\bm{q}_{\mathcal{H}}\right)
			\\=&
			16{\rm Re}\left\{\nabla_{\bm{p}^*} f\left(\bm{p}\right)^{\mathsf{H}} \nabla_{\bm{q}^*} f\left(\bm{q}\right)\right\};
		\end{aligned}\label{tidupq}\\
		&(c)~ \lVert{\nabla_{\mathcal{R}} f\left(\bm{p}_{\mathcal{R}}\right)}\rVert_2
		=
		2 \lVert{\nabla_{\mathcal{H}}f\left(\bm{p}_{\mathcal{H}}\right)}\rVert_2
		=
		4\lVert{\nabla_{\bm{p}} f\left(\bm{p}\right)}\rVert_2. \label{tidufs}
	\end{align}
\end{proposition}
\begin{IEEEproof} By the relationship of  $\bm{q}$, $\bm{q}_{\mathcal{R}}$, and $\bm{q}_{\mathcal{H}}$, and the relationship of $\nabla_{\bm{q}^*} f$, $\nabla_{\mathcal{R}} f$, and $\nabla_{\mathcal{H}^*} f$, we have
\begin{align}
	&\begin{aligned}
		(a)~\nabla_{\mathcal{R}} f\left(\bm{p}_{\mathcal{R}}\right)^{\mathsf{T}} \bm{q}_{\mathcal{R}}
		&=
		\nabla_{\mathcal{R}} f\left(\bm{p}_{\mathcal{R}}\right)^{\mathsf{H}}\bm{q}_{\mathcal{R}}
		\\&\overset{\eqref{r,h}\eqref{grad}}{=}
		\nabla_{\mathcal{H}^*} f\left(\bm{p}_{\mathcal{H}}\right) ^{\mathsf{H}} {\bm{J}_n} \dfrac{1}{4} \bm{J}_n^{\mathsf{H}} \bm{q}_{\mathcal{H}}\quad
		\\&=
		\nabla_{\mathcal{H}^*} f\left(\bm{p}_{\mathcal{H}}\right)^{\mathsf{H}} \bm{q}_{\mathcal{H}}\\
		&=
		\sum\limits_{\mu \in \left\{1,i,j,k \right\}} \nabla_{\bm{p}^*} f\left(\bm{p}\right)^{\mu \mathsf{H}} \bm{q}^{\mu} 
		\\&\overset{\eqref{xzp}}{=}
		\sum\limits_{\mu \in \left\{1,i,j,k \right\}} \left(\nabla_{\bm{p}^*} f\left(\bm{p}\right)^\mathsf{H} \bm{q}\right)^{\mu}
		\\&\overset{\eqref{r,h}}{=}
		4{\rm Re}\left\{\nabla_{\bm{p}^*} f\left(\bm{p}\right)^\mathsf{H} \bm{q}\right\};
	\end{aligned}
\end{align}
\begin{align}
	&\begin{aligned}
		(b)~\nabla_{\mathcal{R}} f\left(\bm{p}_{\mathcal{R}}\right)&^{\mathsf{T}} \nabla_{\mathcal{R}} f\left(\bm{q}_{\mathcal{R}}\right)
		=
		\nabla_{\mathcal{R}} f\left(\bm{p}_{\mathcal{R}}\right)^{\mathsf{H}}\nabla_{\mathcal{R}} f\left(\bm{q}_{\mathcal{R}}\right)\\
		&\overset{\eqref{grad}}{=}
		\nabla_{\mathcal{H}^*} f\left(\bm{p}_{\mathcal{H}}\right) ^{\mathsf{H}} {\bm{J}_n} \bm{J}_n^{\mathsf{H}} \nabla_{\mathcal{H}^*} f\left(\bm{q}_{\mathcal{H}}\right)
		\\&=
		4\nabla_{\mathcal{H}^*} f\left(\bm{p}_{\mathcal{H}}\right)^{\mathsf{H}} \nabla_{\mathcal{H}^*} f\left(\bm{q}_{\mathcal{H}}\right)\\
		&=
		4\sum\limits_{\mu \in \left\{1,i,j,k \right\}} \nabla_{\bm{p}^*} f\left(\bm{p}\right)^{\mu \mathsf{H}} \nabla_{\bm{q}^*} f\left(\bm{q}\right)^{\mu} 
		\\&\overset{\eqref{xzp}}{=}
		4\sum\limits_{\mu \in \left\{1,i,j,k \right\}} \left(\nabla_{\bm{p}^*} f\left(\bm{p}\right)^\mathsf{H} \nabla_{\bm{q}^*} f\left(\bm{q}\right)\right)^{\mu}\\
		&\overset{\eqref{r,h}}{=}
		16{\rm Re}\left\{\nabla_{\bm{p}^*} f\left(\bm{p}\right)^{\mathsf{H}} \nabla_{\bm{q}^*} f\left(\bm{q}\right)\right\};
	\end{aligned}\\
	&\begin{aligned}
		(c)~ \text{Let } \bm{q}=\bm{p} \text{ in } (\ref{tidupq}).\nonumber
	\end{aligned}
\end{align}
This completes the proof.
\end{IEEEproof}
\section{Discriminant Theorems for Convex Quaternion Functions}
The objective of this section is to introduce five discriminant criteria for convex quaternion functions, including the first-order characterization and the second-order characterization. An example is presented to illustrate how these criteria can be applied in practice. 
\subsection{Convex Set and Convex Quaternion Function} 
We begin by introducing the fundamental concepts, such as convex set and convex function \cite{boyd_vandenberghe_2004,nesterov2018lectures}.
%\begin{definition}
%	[Appropriate function]
%	Given the generalized real-valued function $f$ and the nonempty set $\mathcal{X}\subset\bm{dom} f$,  
%	the function $ f $ is called appropriate for the set $ \mathcal{X} $, if $ \exists \bm{x} \in \mathcal{X} $, $ f\left(\bm{x}\right)<+\infty $,  and $\forall \bm{x} \in \mathcal{X}$, $ f\left(\bm{x}\right)>-\infty $.
%	%(????)
%	%???????? $f$ ????? $\mathcal{X}\subset\bm{dom} f$ . ???? $ \exists \bm{x} \in \mathcal{X} $ ?? $ f\left(\bm{x}\right)<+\infty $, ?????? $ \bm{x} \in \mathcal{X}$, ?? $ f\left(\bm{x}\right)>-\infty $, ???? ? $ f $ ???? $ \mathcal{X} $ ????.
%\end{definition}
\begin{definition}
	[Convex set]
	The set $\mathcal{C}$ is called convex, if $\forall \bm{x}, \bm{y} \in \mathcal{C}$, $\forall 0 \leqslant \theta \leqslant 1$, 
	%	(??)
	%	?????? $\mathcal{C}$ ?????????? $\mathcal{C}$ ?, ?? $\mathcal{C}$ ???, ?
	$
		 \theta \bm{x}+\left(1-\theta\right) \bm{y} \in \mathcal{C}.
	$ 
	The set $ \mathcal{C} $ can be a subset of $\mathbb{H}^n$, $\mathcal{R}$ or $\mathcal{H}$.
	%	??, $ \mathcal{C} $???$\mathbb{H}^n$???, ????$\mathcal{R}$?$\mathcal{H}$???.
\end{definition}

\begin{definition}
	[Convex function]\label{tf}
	A function $f$ is said to be convex, if $ \bm{dom} f $ is convex, and $\forall \bm{x}, \bm{y} \in \bm{dom} f$, $0 \leqslant \theta \leqslant 1$,
	%		(???)
	%	??? $f$ ?????, ?? $ \bm{dom} f $ ???, ?
	\begin{equation}\label{tu}
		f\big(\theta \bm{x}+\left(1-\theta\right) \bm{y}\big) \leqslant \theta f\left(\bm{x}\right)+\left(1-\theta\right) f\left(\bm{y}\right).
	\end{equation}
The range of the function $f$ is $\mathbb{R}$, and the definition field $\bm{dom} f$ can be a subset of $\mathbb{H}^n$, $\mathcal{R}$ or $\mathcal{H}$.
	%??? $ \bm{x}, \bm{y} \in \bm{dom} f$, $0 \leqslant \theta \leqslant 1 $ ???, ?? $ f $ ????.
	%Moreover, $f$ is said to be strictly convex, if $ \bm{dom} f $ is convex, and $\forall \bm{x}, \bm{y} \in \bm{dom} f, \bm{x}\ne \bm{y}$, $0 < \theta < 1$,
	%	(?????)
	%	??? $f$ ?????, ?? $ \bm{dom} f $ ???, ?
%	\begin{equation}
%		f\big(\theta \bm{x}+\left(1-\theta\right) \bm{y}\big) < \theta f\left(\bm{x}\right)+\left(1-\theta\right) f\left(\bm{y}\right).
%	\end{equation}
\end{definition}
\begin{example}\label{aq=b}
	Consider a quaternion matrix, $\bm{A} \in \mathbb{H}^{m \times n}$, and a quaternion vector, $\bm{b} \in \mathbb{H}^{m}$, then the set $\mathcal{D} \triangleq \left\{\bm{q} \in \mathbb{H}^{n} \mid  \bm{A} \bm{q} = \bm{b}\right\}$ is convex.
\end{example}
\begin{IEEEproof} 
%$\mathcal{D} \triangleq \left\{\bm{q} \in \mathbb{H}^{n} \mid  \bm{A} \bm{q} = \bm{b}\right\}$, 
$\forall \bm{p}, \bm{q} \in \mathcal{D}$, $\bm{A} \bm{p} = \bm{b}$, $\bm{A} \bm{q} = \bm{b}$, 
$\forall 0 \leqslant \theta \leqslant 1$,
\begin{equation}
	\bm{A} \left(\theta \bm{p} + (1 - \theta) \bm{q}\right)
	=
	\theta \bm{A} \bm{p} + (1 - \theta) \bm{A} \bm{q}
	=
	\theta \bm{b} + (1 - \theta) \bm{b}
	=
	\bm{b}.
\end{equation}
Therefore, $\theta \bm{p} + (1 - \theta) \bm{q} \in \mathcal{D}$, that is the set $\mathcal{D}$ is convex.
\end{IEEEproof}
\begin{example}\label{fq}
	If the quaternion function $f(\bm{q})$ is convex, then the set $\mathcal{E} \triangleq \left\{\bm{q} \in \mathbb{H}^{n} \mid  f(\bm{q}) \leqslant 0\right\}$ is also convex.
\end{example}
\begin{IEEEproof} 
%$\mathcal{E} \triangleq \left\{\bm{q} \in \mathbb{H}^{n} \mid  f(\bm{q}) \leqslant 0\right\}$, 
$\forall \bm{p}, \bm{q} \in \mathcal{E}$, $f(\bm{p}) \leqslant 0$, $f(\bm{q}) \leqslant 0$. Since $f(\bm{q})$ is convex, 
$\forall 0 \leqslant \theta \leqslant 1$,
\begin{equation}
	f \big(\theta \bm{p} + (1 - \theta) \bm{q}\big)
	\leqslant
	\theta f(\bm{p}) + (1 - \theta) f(\bm{q})
	\leqslant
	0.
\end{equation}
Therefore, $\theta \bm{p} + (1 - \theta) \bm{q} \in \mathcal{E}$, that is the set $\mathcal{E}$ is convex.
\end{IEEEproof}
%\begin{definition}
%	[Strictly convex function]
%	The appropriate function $f$ is said to be strictly convex, if $ \bm{dom} f $ is convex, and $\forall \bm{x}, \bm{y} \in \bm{dom} f, \bm{x}\ne \bm{y}$, $0 < \theta < 1$,
%	%	(?????)
%	%	??? $f$ ?????, ?? $ \bm{dom} f $ ???, ?
%	\begin{equation}
%		f\big(\theta \bm{x}+\left(1-\theta\right) \bm{y}\big) < \theta f\left(\bm{x}\right)+\left(1-\theta\right) f\left(\bm{y}\right).
%	\end{equation}
%	%	??? $ \bm{x}, \bm{y} \in \bm{dom} f$, $\bm{x}\ne \bm{y}$, $0 < \theta < 1 $ ???, ?? $ f $ ??????.
%\end{definition}
\subsection{First-order Characterization of Discriminant Theorems for Convex Quaternion Functions}
We shall now introduce four discriminant theorems for convex quaternion functions, including the first-order characterization and gradient monotonicity.
\begin{theorem}\label{c,cr,ch}
	Consider the three sets $\mathcal{C} \subset \mathbb{H}^n$, \\$
	\mathcal{C}_{\mathcal{R}} \triangleq  \left\{\bm{q}_{\mathcal{R}}=\left(\bm{q}_a^{\mathsf{T}},\bm{q}_b^{\mathsf{T}},\bm{q}_c^{\mathsf{T}},\bm{q}_d^{\mathsf{T}}\right)^{\mathsf{T}} \mid \bm{q} \in\mathcal{C}\right\}\subset \mathcal{R} = \mathbb{R}^{4n}$, \\
	$\mathcal{C}_{\mathcal{H}} 
	\triangleq 
	\left\{\bm{q}_{\mathcal{H}}=\left(\bm{q}^{\mathsf{T}},\bm{q}^{i\mathsf{T}},\bm{q}^{j\mathsf{T}},\bm{q}^{k\mathsf{T}}\right)^{\mathsf{T}} \mid  \bm{q}\in\mathcal{C}\right\}
	\subset \mathcal{H} \subset \mathbb{H}^{4n}$. Then,  
	$
		\mathcal{C} \text{ is convex} \Leftrightarrow \mathcal{C}_{\mathcal{R}} \text{ is convex} \Leftrightarrow \mathcal{C}_{\mathcal{H}} \text{ is convex} .
	$
	\begin{IEEEproof} 
	Using the definition of $\mathcal{C}$, $\mathcal{C}_{\mathcal{R}}$,  $\mathcal{C}_{\mathcal{H}}$, and that of convex set, the proof following.
	\end{IEEEproof}
\end{theorem}

A straightforward method to discriminate the convexity of a quaternion function is to confine it to a line segment and determine whether the resulting one-dimensional function is convex, as in the following theorem.
\begin{theorem}\label{zx}
	The quaternion function $f(\bm{q}):\mathcal{C}\subset \mathbb{H}^n \to \mathbb{R}$ is convex if and only if (shortened to iff)  $\forall \bm{q} \in \mathcal{C}$, $\bm{v} \in \mathbb{H}^{n}$, $ g:\mathcal{S} \rightarrow \mathbb{R}$, 
	\begin{equation}
		g(t) = f(\bm{q}+t\bm{v})
	\end{equation}
	is convex, 
	where $\mathcal{S} \triangleq \left\{t \in \mathbb{R}\mid  \bm{q} + t\bm{v} \in \mathcal{C}\right\} \subset \mathbb{R}$. 
%	Moreover, $f(\bm{q})$ is strictly convex, iff $\forall \bm{q} \in \mathcal{C}$, $\bm{v} \in \mathbb{H}^{n}$, $g(t)$ is strictly convex.
\end{theorem}
\begin{IEEEproof} 
The proof follows the same steps as its counterpart in the real field \cite{boyd_vandenberghe_2004,nesterov2018lectures}.
\end{IEEEproof}

%\begin{corollary}
%	The quaternion function $f(\bm{q}):\mathcal{C} \subset \mathbb{H}^n \to \mathbb{R}$ is strictly convex, iff $\forall \bm{q} \in \mathcal{C}$, $\bm{v} \in \mathbb{H}^{n}$, $\mathcal{S} \triangleq \left\{t \in \mathbb{R}\mid  \bm{q} + t\bm{v} \in \mathcal{C}\right\} \subset \mathbb{R}$, $ g:\mathcal{S} \rightarrow \mathbb{R}$, where 
%	\begin{equation}
%		g(t) = f(\bm{q}+t\bm{v})
%	\end{equation}
%	is strictly convex.
%\end{corollary}
%\begin{IEEEproof} 
%The proof is similar to that of Theorem $\ref{zx}$.
%\end{IEEEproof}
For real-differentiable quaternion functions, we can also use their gradient information to discriminate their convexity, as stated in the following theorem.
\begin{theorem}
	[First-order characterization \cite{flamant2021general}]\label{yijie}
	Consider a convex set $\mathcal{C}\subset \mathbb{H}^n$ and a real-differentiable quaternion function $f(\bm{q}):\mathcal{C} \to \mathbb{R}$. Then $f(\bm{q})$ is convex iff 
	$\forall \bm{p}, \bm{q} \in \mathcal{C}$,
	\begin{equation}\label{yj}
		f\left(\bm{q}\right) \geqslant f\left(\bm{p}\right)+4{\rm Re}\left\{\nabla_{\bm{p}^*}f\left(\bm{p}\right)^{\mathsf{H}}\left(\bm{q}-\bm{p}\right)\right\},
	\end{equation}
%	Moreover, $f(\bm{q})$ is strictly convex iff 
%	$\forall \bm{p}, \bm{q} \in \mathcal{C}$, $\bm{p} \ne \bm{q}$,
%	\begin{equation}
%		f\left(\bm{q}\right) > f\left(\bm{p}\right) + 4{\rm Re}\left\{\nabla_{\bm{p}^*} f\left(\bm{p}\right)^{\mathsf{H}} \left(\bm{q}-\bm{p}\right)\right\},
%	\end{equation}
	where $\nabla_{\bm{p}^*}f\left(\bm{p}\right)$ is defined in $(\ref{getd})$.
\end{theorem}	

Another commonly used first-order characterization is gradient monotonicity, as shown below.
\begin{theorem}[Gradient monotonicity]\label{dandiao}
	Consider a convex set $\mathcal{C}\subset \mathbb{H}^n$ and a real-differentiable quaternion function $f(\bm{q}):\mathcal{C} \to \mathbb{R}$.
	Then, $f(\bm{q})$ is convex iff $\forall \bm{p}, \bm{q} \in \mathcal{C}$,
	%	(?????) ?????$f(\bm{q}):\mathcal{C}\subset \mathbb{H}^n \to \mathbb{R}$??????, ? $ f(\bm{q}) $ ???????? $\mathcal{C}$ ????
	\begin{equation}\label{dandiaogs}
		{\rm Re}\left\{\big(\nabla_{\bm{p}^*} f\left(\bm{p}\right)-\nabla_{\bm{q}^*} f\left(\bm{q}\right)\big)^{\mathsf{H}}\left(\bm{p}-\bm{q}\right)\right\} \geqslant 0,
	\end{equation}
%	Moreover, $f(\bm{q})$ is strictly convex iff 
%	$\forall \bm{p}, \bm{q} \in \mathcal{C}$, $\bm{p} \ne \bm{q}$,
%	\begin{equation}
%		{\rm Re}\left\{\big(\nabla_{\bm{p}^*} f\left(\bm{p}\right)-\nabla_{\bm{q}^*} f\left(\bm{q}\right)\big)^{\mathsf{H}}\left(\bm{p}-\bm{q}\right)\right\} > 0,
%	\end{equation}
	where $\nabla_{\bm{p}^*}f\left(\bm{p}\right)$ is defined in $(\ref{getd})$.
	
	%??, $\nabla_{\bm{p}^*}f\left(\bm{p}\right)$?$f\left(\bm{p}\right)$?????, ?$(\ref{getd})$???.
\end{theorem}
\begin{IEEEproof} 
From Theorem $\ref{c,cr,ch}$, $\mathcal{C}$ is convex iff $\mathcal{C}_{\mathcal{R}}$ is convex. 
We already know  \cite{boyd_vandenberghe_2004,nesterov2018lectures} that for a differentiable real function, 
$f(\bm{q}_{\mathcal{R}})$ is convex iff $\forall \bm{p}_{\mathcal{R}}, \bm{q}_{\mathcal{R}} \in \mathcal{C}_{\mathcal{R}}$,
\begin{equation}\label{dd}
	\big(\nabla_{\mathcal{R}}f\left(\bm{p}_{\mathcal{R}}\right)-\nabla_{\mathcal{R}}f\left(\bm{q}_{\mathcal{R}}\right)\big)^{\mathsf{T}}\left(\bm{p}_{\mathcal{R}}-\bm{q}_{\mathcal{R}}\right) \geqslant 0,
\end{equation}
where the set $\mathcal{C}_{\mathcal{R}} \subset \mathcal{R}$ is convex. 
Hence from  $\left(\ref{tiduq}\right)$, we have
%???$\ref{c,cr,ch}$?, $\mathcal{C}_{\mathcal{R}}$??????$\mathcal{C}$???.?$\left(\ref{tiduq}\right)$?
\begin{align}\label{tdp}
	\begin{aligned}
		&\big(\nabla_{\mathcal{R}}f\left(\bm{p}_{\mathcal{R}}\right)-\nabla_{\mathcal{R}}f\left(\bm{q}_{\mathcal{R}}\right)\big)^{\mathsf{T}}\left(\bm{p}_{\mathcal{R}}-\bm{q}_{\mathcal{R}}\right)
		\\=&
		4{\rm Re}\left\{\big(\nabla_{\bm{p}^*} f\left(\bm{p}\right)-\nabla_{\bm{q}^*} f\left(\bm{q}\right)\big)^{\mathsf{H}}\left(\bm{p}-\bm{q}\right)\right\}.
	\end{aligned}
\end{align}
Upon substituting $(\ref{tdp})$ into $(\ref{dd})$, the proof follows. 
%In a similar manner, we can prove the equivalence of strictly convex functions.
\end{IEEEproof}
%\begin{corollary}[Gradient monotonicity]\label{tlqtu}
%	Consider a convex set $\mathcal{C}\subset \mathbb{H}^n$ and a real-differentiable quaternion function $f(\bm{q}):\mathcal{C} \to \mathbb{R}$. Then, $f(\bm{q})$ is strictly convex iff 
%	$\forall \bm{p}, \bm{q} \in \mathcal{C}$, $\bm{p} \ne \bm{q}$,
%	\begin{equation}
%		{\rm Re}\left\{\big(\nabla_{\bm{p}^*} f\left(\bm{p}\right)-\nabla_{\bm{q}^*} f\left(\bm{q}\right)\big)^{\mathsf{H}}\left(\bm{p}-\bm{q}\right)\right\} > 0,
%	\end{equation}
%	 in equation $(\ref{getd})$.
%	%??, $\nabla_{\bm{p}^*}f\left(\bm{p}\right)$?$f\left(\bm{p}\right)$?????, ?$(\ref{getd})$???.
%\end{corollary}
%\begin{IEEEproof} 
%The proof is similar to that of Theorem $\ref{dandiao}$.
%\end{IEEEproof}

In addition, we can also use the epigraph to discriminate the convexity of $f\left(\bm{q}\right)$, as shown below.
\begin{definition}[Epigraph]
	For the quaternion generalized real-valued function $f\left(\bm{q}\right):\mathbb{H}^n \to \mathbb{R}\cup \left\{\pm \infty \right\}$, the set 
	%???????????$f\left(\bm{q}\right):\mathbb{H}^n \to \mathbb{R}\cup \left\{\pm \infty \right\}$, 
	\begin{equation}
		\bm{epi} f = \left\{(\bm{q},\, t) \in \mathbb{H}^{n+1} \mid f(\bm{q}) \leqslant t, t \in \mathbb{R}\right\}
	\end{equation}
	is called the epigraph of $ f(\bm{q}) $.
	%?? $ f(\bm{q}) $ ????.
\end{definition}
\begin{theorem}\label{sft}
	The quaternion generalized real-valued function $f\left(\bm{q}\right):\mathcal{C} \subset \mathbb{H}^n \to \mathbb{R}\cup \left\{\pm \infty \right\}$ is convex, iff $\bm{epi} f$ is a convex set.
	%	?????????$f\left(\bm{q}\right):\mathcal{C} \subset \mathbb{H}^n \to \mathbb{R}\cup \left\{\pm \infty \right\}$????????????$\bm{epi} f$ ???. 
\end{theorem}
\begin{IEEEproof} 
	The proof follows the same steps as its counterpart in the real field \cite{boyd_vandenberghe_2004,nesterov2018lectures}.
\end{IEEEproof}
\subsection{Second-order Characterization of Discriminant Theorems for Convex Quaternion Functions}
Before introducing the second-order characterization of convex quaternion functions, we first need to define positive definite quaternion matrices.
\begin{definition}[Positive definite matrix]\label{zdjz}
	The matrix $\bm{A}\in\mathbb{H}^{n\times n}$ is called positive definite, if
	%	(????)
	%	?$\bm{A}\in\mathbb{H}^{n\times n}$, ????$\bm{x} \in \mathbb{H}^n$?$\bm{x} \ne \bm{0} $, ??
	\begin{equation}
		{\rm Re}\left\{\bm{x}^{\mathsf{H}} \bm{A} \bm{x}\right\} >0,\quad  \forall \bm{x} \in \mathbb{H}^n,\, \bm{x} \ne \bm{0},
	\end{equation}
	and is denoted by $\bm{A} \succ \bm{O}$, where $\bm{O}$ is the $n \times n$ zero matrix. Similarly, the matrix $\bm{A}\in\mathbb{H}^{n\times n}$ is called positive semi-definite, if 
	\begin{equation}
		{\rm Re}\left\{\bm{x}^{\mathsf{H}} \bm{A} \bm{x}\right\} \geqslant 0,\quad \forall \bm{x} \in \mathbb{H}^n,\, \bm{x} \ne \bm{0},
	\end{equation}
	and is denoted by $\bm{A} \succeq \bm{O}$.
	%	??$\bm{A}$?????, ??$\bm{A} \succ \bm{O}$, ??$\bm{O}$????.\quad ?${\rm Re}\left\{\bm{x}^{\mathsf{H}} \bm{A} \bm{x}\right\}\geqslant 0$, ??$\bm{A}$??????, ??$\bm{A} \succeq \bm{O}$.
\end{definition}
\begin{theorem}\label{dczd}
	If the matrix $\bm{A}\in\mathbb{H}^{n\times n}$ satisfies $\bm{A}^{\mathsf{H}}=\bm{A}$, then $\bm{A}$ is positive definite iff 
	%	?$\bm{A}\in\mathbb{H}^{n\times n}$, ?$\bm{A}^{\mathsf{H}}=\bm{A}$, ?$\bm{A}$????????
	\begin{equation}
		\bm{x}^{\mathsf{H}} \bm{A} \bm{x}>0,\quad \forall \bm{x} \in \mathbb{H}^n,\, \bm{x} \ne \bm{0}.
	\end{equation}
	Similarly, the matrix $\bm{A}$ is positive semi-definite iff 
	%$\bm{A}$?????????
	\begin{equation}
		\bm{x}^{\mathsf{H}} \bm{A} \bm{x}\geqslant 0,\quad \forall \bm{x} \in \mathbb{H}^n,\,\bm{x} \ne \bm{0}.
	\end{equation}
\end{theorem}
\begin{IEEEproof} 
This is straightforward to prove, by applying Definition $\ref{zdjz}$.	
%$\forall \bm{x} \in \mathbb{H}^n$, $\bm{x} \ne \bm{0}$, since $\bm{A}^{\mathsf{H}}=\bm{A}$, so that
%\begin{equation}
%	\bm{x}^{\mathsf{H}} \bm{A} \bm{x}
%	=
%	\bm{x}^{\mathsf{H}} \bm{A}^{\mathsf{H}} \bm{x}
%	=
%	\left(\bm{x}^{\mathsf{H}} \bm{A} \bm{x}\right)^{\mathsf{H}}, 
%\end{equation}
%that is $\bm{x}^{\mathsf{H}} \bm{A} \bm{x} \in \mathbb{R}$. Therefore,
%\begin{equation}
%	\bm{x}^{\mathsf{H}} \bm{A} \bm{x}
%	=
%	{\rm Re}\left\{\bm{x}^{\mathsf{H}} \bm{A} \bm{x}\right\},\quad \forall \bm{x} \in \mathbb{H}^n,\, \bm{x} \ne \bm{0}.
%\end{equation}
%This completes the proof.
\end{IEEEproof}

If the quaternion function $f(\bm{q})$ is second-order continuous real-differentiable, we can use the Hessian matrix to discriminate its convexity, as shown below.
\begin{theorem}[Second-order characterization]\label{er}
	Consider a convex set $\mathcal{C}\subset \mathbb{H}^n$ and a second-order continuous real-differentiable quaternion function $f(\bm{q}):\mathcal{C} \to \mathbb{R}$. Then $f(\bm{q})$ is convex iff 
	%(????)?????$f(\bm{q}):\mathcal{C} \to \mathbb{R}$??????????, ??$\mathcal{C}\subset \mathbb{H}^n$???, ? $ f $ ????????
	\begin{equation}
		\boldsymbol{H}_{\mathcal{H}\mathcal{H}^*} \succeq \bm{O}.
	\end{equation}
%	In addition, $f(\bm{q})$ is strictly convex if 
%	%$ f $ ??????????
%	\begin{equation}
%		\boldsymbol{H}_{\mathcal{H}\mathcal{H}^*} \succ \bm{O},
%	\end{equation}
	where $\boldsymbol{H}_{\mathcal{H}\mathcal{H}^*}$ is defined in $(\ref{HHH})$.
\end{theorem}
\begin{IEEEproof} 
Applying Theorem $\ref{c,cr,ch}$, the set $\mathcal{C}$ is convex iff the set $\mathcal{C}_{\mathcal{R}}$ is convex. 
%???$\ref{c,cr,ch}$?, $\mathcal{C}$??????$\mathcal{C}_{\mathcal{R}}$???. ?$(\ref{iff})$?, $f\left(\bm{q}\right)$ ???$ f\left(\bm{q}_{\mathcal{R}}\right)$.
We already know \cite{boyd_vandenberghe_2004,nesterov2018lectures} that for a second-order continous differentiable function, 
$f(\bm{q}_{\mathcal{R}})$ is convex iff
\begin{equation}\label{hession0}
	\boldsymbol{H}_{\mathcal{R}\mathcal{R}} \succeq  \bm{O}, \quad \forall \bm{q}_{\mathcal{R}} \in \mathcal{C}_{\mathcal{R}},
\end{equation}
where the set $\mathcal{C}_{\mathcal{R}} \subset \mathcal{R}$ is convex. 
By Corollary $\ref{hermite}$,  $\boldsymbol{H}_{\mathcal{H}\mathcal{H}^*}$ is a Hermite matrix. Then $\forall \bm{x}_{\mathcal{H}} \in \mathcal{H}$, $\bm{x}_{\mathcal{H}} \ne \bm{0}$, we have
%	?$\left(\ref{hassion}\right)$?$\boldsymbol{H}_{\mathcal{H}\mathcal{H}^*}=\dfrac{1}{16}\bm{J}_n\boldsymbol{H}_{\mathcal{R}\mathcal{R}}\bm{J}_n^{\mathsf{H}}$, ?$\boldsymbol{H}_{\mathcal{H}\mathcal{H}^*}$?Hermite??, ?$\forall \bm{x}_{\mathcal{H}} \in \mathcal{H}$?$\bm{x}_{\mathcal{H}}  \ne \bm{0}$, 
\begin{align}
	\begin{aligned}
		\bm{x}_{\mathcal{H}} ^{\mathsf{H}} \boldsymbol{H}_{\mathcal{H}\mathcal{H}^*} \bm{x}_{\mathcal{H}} 
		\overset{\eqref{hassion}}{=} 
		&\dfrac{1}{16} \bm{x}_{\mathcal{H}} ^{\mathsf{H}} \bm{J}_n \boldsymbol{H}_{\mathcal{R}\mathcal{R}} \bm{J}_n^{\mathsf{H}} \bm{x}_{\mathcal{H}} 
		\\=&
		\dfrac{1}{16} \left(\bm{J}_n^{\mathsf{H}} \bm{x}_{\mathcal{H}} \right)^{\mathsf{H}} \boldsymbol{H}_{\mathcal{R}\mathcal{R}} \left(\bm{J}_n^{\mathsf{H}} \bm{x}_{\mathcal{H}} \right)
		\\\overset{\eqref{r,h}}{=}&
		\bm{x}_{\mathcal{R}}^{\mathsf{H}} \boldsymbol{H}_{\mathcal{R}\mathcal{R}} \bm{x}_{\mathcal{R}}.
	\end{aligned}
\end{align}
Therefore,
\begin{equation}
	\boldsymbol{H}_{\mathcal{H}\mathcal{H}^*}\succeq  \bm{O} \quad \Leftrightarrow \quad \boldsymbol{H}_{\mathcal{R}\mathcal{R}} \succeq  \bm{O},
\end{equation}
which concludes the proof. 
%When the function $f(\bm{q})$ is strictly convex, 
%the proof is similar to the above.
%$ f $ ???????, ????.
\end{IEEEproof}
\begin{corollary}\label{ertuilun}
	Consider a convex set $\mathcal{C}\subset \mathbb{H}^n$ and a second-order continuous real-differentiable quaternion function $f(\bm{q}):\mathcal{C} \to \mathbb{R}$. Then, the following three propositions are equivalent:
	%	?????$f(\bm{q}):\mathcal{C} \to \mathbb{R}$??????????, ??$\mathcal{C}\subset \mathbb{H}^n$???, ?????????
	
	$\mathit{(a)}$~ $ f(\bm{q}) $ is convex;
	
	$\mathit{(b)}$~ $\boldsymbol{H}_{\mathcal{H}\mathcal{H}^*}\succeq  \bm{O}$;
	
	$\mathit{(c)}$~ $\sum\limits_{\nu \in \left\{1,i,j,k \right\}} {\rm Re}\left\{\bm{x}^{\mathsf{H}} \boldsymbol{H}_{\bm{q}^\nu \bm{q}^{*}} \bm{x}^{\nu}\right\}
	\geqslant 0
	$,~ $\forall \bm{x} \in \mathbb{H}^n$, $\bm{x} \ne \bf{0} .$
\end{corollary}
\begin{IEEEproof} 
From Theorem
$\ref{er}$, $(a)$ is equivalent to $(b)$, so we only need to prove that $(b)$ is equivalent to $(c)$. 
From Corollary $\ref{hermite}$, we know that $\boldsymbol{H}_{\mathcal{H}\mathcal{H}^*}$ is a Hermite matrix. Then $\forall \bm{x}_{\mathcal{H}} \in \mathcal{H}$, $\bm{x}_{\mathcal{H}} \ne \bm{0}$, we have
%??$(\ref{Hermite})$,$\boldsymbol{H}_{\mathcal{H}\mathcal{H}^*}$?Hermite????, ?
\begin{align}\label{xhx}
	&\bm{x}_{\mathcal{H}} ^{\mathsf{H}} \boldsymbol{H}_{\mathcal{H}\mathcal{H}^*} \bm{x}_{\mathcal{H}} \nonumber\\
	=&
	\begin{pmatrix}
		\bm{x}	\\
		\bm{x}^i	\\
		\bm{x}^j	\\
		\bm{x}^k	
	\end{pmatrix}^{\mathsf{H}}
	%	\begin{pmatrix}
		%		\bm{x}^{\mathsf{H}},\bm{x}^{i\mathsf{H}},\bm{x}^{j\mathsf{H}},\bm{x}^{k\mathsf{H}}
		%	\end{pmatrix}\\&
	\left(\begin{array}{llll}
		\boldsymbol{H}_{\bm{q} \bm{q}^{*}} & \boldsymbol{H}_{\bm{q}^{i} \bm{q}^{*}} & \boldsymbol{H}_{\bm{q}^{j} \bm{q}^{*}} & \boldsymbol{H}_{\bm{q}^{k} \bm{q}^{*}} \\
		\boldsymbol{H}_{\bm{q} \bm{q}^{i *}} & \boldsymbol{H}_{\bm{q}^{i} \bm{q}^{i *}} & \boldsymbol{H}_{\bm{q}^{j} \bm{q}^{i *}} & \boldsymbol{H}_{\bm{q}^{k} \bm{q}^{i *}} \\
		\boldsymbol{H}_{\bm{q} \bm{q}^{j *}} & \boldsymbol{H}_{\bm{q}^{i} \bm{q}^{j *}} & \boldsymbol{H}_{\bm{q}^{j} \bm{q}^{j *}} & \boldsymbol{H}_{\bm{q}^{k} \bm{q}^{j *}} \\
		\boldsymbol{H}_{\bm{q} \bm{q}^{k *}} & \boldsymbol{H}_{\bm{q}^{i} \bm{q}^{k *}} & \boldsymbol{H}_{\bm{q}^{j} \bm{q}^{k *}} & \boldsymbol{H}_{\bm{q}^{k} \bm{q}^{k *}}
	\end{array}\right)
	\begin{pmatrix}
		\bm{x}	\\
		\bm{x}^i	\\
		\bm{x}^j	\\
		\bm{x}^k	
	\end{pmatrix}\nonumber
	\\=&
	%		\sum_{\nu,\mu \in \left\{1,i,j,k \right\}} \bm{x}^{\mu \mathsf{H}} \dfrac{\partial}{\partial \bm{q}^\nu}\left(\dfrac{\partial f}{\partial \bm{q}^{\mu *}}\right)^{\mathsf{T}} \bm{x}^{\nu}
	%		=
	\sum_{\mu,\nu \in \left\{1,i,j,k \right\}} \bm{x}^{\mu \mathsf{H}} \boldsymbol{H}_{\bm{q}^\nu \bm{q}^{\mu*}} \bm{x}^{\nu}\nonumber
	\\\overset{\eqref{r,h}}{=}&
	%		4\sum_{\mu \in \left\{1,i,j,k \right\}} {\rm Re}\left(\bm{x}^{\mu \mathsf{H}} \boldsymbol{H}_{\bm{q} \bm{q}^{\mu*}} \bm{x}\right)
	%		=
	4\sum_{\nu \in \left\{1,i,j,k \right\}} {\rm Re}\left\{\bm{x}^{\mathsf{H}} \boldsymbol{H}_{\bm{q}^\nu \bm{q}^{*}} \bm{x}^{\nu}\right\}.
\end{align}
%	\left(??
%	\begin{equation}
	%		\bm{x}_{\mathcal{H}} ^{\mathsf{H}} \boldsymbol{H}_{\mathcal{H}\mathcal{H}^*} \bm{x}_{\mathcal{H}} 
	%		= 
	%		16\bm{x}^{\mathsf{H}} \boldsymbol{H}_{\bm{q}\bm{q}^*} \bm{x}
	%	\end{equation}\right)
Therefore,
\begin{align}
	\begin{aligned}
		\sum_{\nu \in \left\{1,i,j,k \right\}} {\rm Re}&\left\{\bm{x}^{\mathsf{H}} \boldsymbol{H}_{\bm{q}^\nu \bm{q}^{*}} \bm{x}^{\nu}\right\}\geqslant 0, \quad \forall \bm{x} \in \mathbb{H}^n, \,  \bm{x} \ne \bm{0}\\ &\quad \Leftrightarrow \quad \boldsymbol{H}_{\mathcal{H}\mathcal{H}^*}\succeq  \bm{O}.
	\end{aligned}
\end{align}
This completes the proof.
\end{IEEEproof}
%\begin{corollary}
%	Consider a convex set $\mathcal{C}\subset \mathbb{H}^n$ and a second-order continuous real-differentiable quaternion function $f(\bm{q}):\mathcal{C} \to \mathbb{R}$. Then, the following three propositions are equivalent:
%	%?????$f(\bm{q}):\mathcal{C} \to \mathbb{R}$??????????, ??$\mathcal{C}\subset \mathbb{H}^n$???, ?????????
%	
%	$\mathit{(a)}$~ $f(\bm{q}) $ is strictly convex;
%	
%	$\mathit{(b)}$~ $\boldsymbol{H}_{\mathcal{H}\mathcal{H}^*}\succ  \bm{O}$;
%	
%	$\mathit{(c)}$~ $\sum\limits_{\nu \in \left\{1,i,j,k \right\}} {\rm Re}\left\{\bm{x}^{\mathsf{H}} \boldsymbol{H}_{\bm{q}^\nu \bm{q}^{*}} \bm{x}^{\nu}\right\}
%	> 0
%	$, $\forall \bm{x} \in \mathcal{C}$, $\bm{x} \ne \bm{0} .$
%\end{corollary}
%\begin{IEEEproof} 
%The proof is similar to that of Corollary $\ref{ertuilun}.$
%\end{IEEEproof}
\begin{lemma}\label{zhuzishi}
	The matrix $\bm{A}\in\mathbb{H}^{n\times n}$ is positive definite (positive semi-definite), iff all principal submatrices of $\bm{A}$ are positive definite (positive semi-definite).
	%?$\bm{A}\in\mathbb{H}^{n\times n}$, ?$\bm{A}$?????(?????)???$\bm{A}$???????????(?????).
\end{lemma}
\begin{IEEEproof} 
The follows the same steps as its counterpart in the real field \cite{horn2012matrix}.
\end{IEEEproof}
%\begin{corollary}\label{zhengding}
%	If the augmented quaternion Hessian matrix $\boldsymbol{H}_{\mathcal{H}\mathcal{H}^*} \in  \mathbb{H}^{4n\times 4n}$ is positive define (positive semi-definite), then the quaternion Hessian matrix $\boldsymbol{H}_{\bm{q}\bm{q}^*}\in \mathbb{H}^{n\times n}$ is positive define (positive semi-definite).
%	%	?$\boldsymbol{H}_{\mathcal{H}\mathcal{H}^*}\in \mathbb{H}^{4n\times 4n}$?????(?????), ?$\boldsymbol{H}_{\bm{q}\bm{q}^*}\in \mathbb{H}^{n\times n}$?????(?????).
%\end{corollary}
%\begin{IEEEproof} 
%This is straightforward to prove from $(\ref{HHH})$ and Lemma $\ref{zhuzishi}$.
%%$\boldsymbol{H}_{\bm{q}\bm{q}^*}\in \mathbb{H}^{n\times n}$?$\boldsymbol{H}_{\mathcal{H}\mathcal{H}^*}\in \mathbb{H}^{4n\times 4n}$????, ????$\ref{zhuzishi}$??.
%\end{IEEEproof}
%????3.1?????????????????????

Applying Lemma $\ref{zhuzishi}$, we can obtain a necessary condition for convex quaternion functions.
\begin{theorem}
	Consider a convex set $\mathcal{C}\subset \mathbb{H}^n$ and a second-order continuous real-differentiable quaternion function $f(\bm{q}):\mathcal{C} \to \mathbb{R}$. If $f(\bm{q})$ is convex, then 
	%the quaternion Hessian matrix $\boldsymbol{H}_{\bm{q}\bm{q}^*}$ is positive semi-definite.
	\begin{equation}
		\boldsymbol{H}_{\bm{q}\bm{q}^*} \succeq \bm{O},
	\end{equation}
	%	In addition, $f(\bm{q})$ is strictly convex if 
	%	%$ f $ ??????????
	%	\begin{equation}
		%		\boldsymbol{H}_{\mathcal{H}\mathcal{H}^*} \succ \bm{O},
		%	\end{equation}
	where $\boldsymbol{H}_{\bm{q}\bm{q}^*}$ is the quaternion Hessian matrix, defined in $(\ref{hqq*})$.
	%?????$f(\bm{q}):\mathcal{C} \to \mathbb{R}$??????????, ??$\mathcal{C}\subset \mathbb{H}^n$???, ? $ f $ ????, ?$\boldsymbol{H}_{\bm{q}\bm{q}^*}$???????? $ f $ ??????, ?$\boldsymbol{H}_{\bm{q}\bm{q}^*}$?????.
\end{theorem}
\begin{IEEEproof} 
	Upon applying Theorem $\ref{er}$, together with the convexity of $f(\bm{q})$, we have $\boldsymbol{H}_{\mathcal{H}\mathcal{H}^*} \succeq \bm{O}$. By $(\ref{HHH})$ and Lemma $\ref{zhuzishi}$, we finally obtain
	$\boldsymbol{H}_{\bm{q}\bm{q}^*} \succeq \bm{O}$.
%	???$\ref{er}$???, $ f $ ???????$\boldsymbol{H}_{\mathcal{H}\mathcal{H}^*}\succeq  \bm{O}$.????$\ref{zhengding}$?, $\boldsymbol{H}_{\bm{q}\bm{q}^*}$??????.\quad $ f $ ???????, ????.
\end{IEEEproof}
\subsection{Examples of Convex Quaternion Function}
We next provide a basic example to determine the convexity of quaternion functions. In this example, we make use of certain GHR derivatives presented in TABLE IV of \cite{Mandic2015The}, which are included in TABLE $\mathrm{\ref{table}}$ here.
\begin{table}[ht]
	\renewcommand\arraystretch{2}
	% The table caption is above the table.
	\vspace{.1 in}
	\caption{Several Derivatives Performed by the GHR Calculus from TABLE IV of \cite{Mandic2015The}, $\forall \boldsymbol{A} \in \mathbb{H}^{n \times n}$, $\forall \boldsymbol{a} \in \mathbb{H}^{n}$, $\forall \boldsymbol{b} \in \mathbb{H}^{n}$, $\alpha \in \mathbb{H}$, $\beta \in \mathbb{H}$.}
	\label{table}
	\begin{center}
		\begin{tabular}{|c|c|c|}
			\hline
			$f(\bm{q})$ or $\bm{f}(\bm{q})$ & $\frac{\partial f}{\partial \bm{q}}$ or $\frac{\partial \bm{f}}{\partial \bm{q}}$ & $\frac{\partial f}{\partial \bm{q}^*}$ or $\frac{\partial \bm{f}}{\partial \bm{q}^*}$\\
			\hline
			$\bm{a}^{\mathsf{T}}\bm{q}\beta$ & $\bm{a}^{\mathsf{T}}{\rm Re}\{\beta\}$ & $-\frac{1}{2}\bm{a}^{\mathsf{T}}\beta^*$\\
			%\hline
			$\alpha \boldsymbol{q}^{\mathsf{H}} \boldsymbol{b}$ & $-\frac{1}{2} \alpha \boldsymbol{b}^{\mathsf{H}}$ & $\alpha {\rm Re}\left\{\boldsymbol{b}^{\mathsf{T}}\right\}$\\
			%\hline
			$\boldsymbol{A} \boldsymbol{q} \beta$ & $\boldsymbol{A}{\rm Re}\left\{\beta\right\}$ & $-\frac{1}{2} \boldsymbol{A} \beta^{*}$\\
			%\hline
			$\boldsymbol{q}^{\mathsf{H}} \boldsymbol{A} \boldsymbol{q}$ & $\boldsymbol{q}^{\mathsf{H}} \boldsymbol{A}-\frac{1}{2}\left(\boldsymbol{A} \boldsymbol{q}\right)^{\mathsf{H}}$ & $-\frac{1}{2} \boldsymbol{q}^{\mathsf{H}} \boldsymbol{A}+{\rm Re}\left\{\left(\boldsymbol{A} \boldsymbol{q}\right)^{\mathsf{T}}\right\}$\\
			\hline
		\end{tabular}
	\end{center}
\end{table}

\begin{example}\label{aq-b}
	If the quaternion function 
	$f\left(\bm{q}\right)=\lVert{\bm{A}\bm{q}-\bm{b}}\rVert_2^2$, $\forall \bm{q} \in \mathbb{H}^n$, $\bm{A}\in \mathbb{H}^{m\times n}$, $\bm{b}\in \mathbb{H}^m$, then $f\left(\bm{q}\right)$ is convex.
\end{example}
\begin{IEEEproof}
(\textbf{First-order characterization criterion})
By the definition of the 2-norm, we have
%, applying Theorem $\ref{yijie}$)
 \begin{align}\label{62}
 	\begin{aligned}
 		f\left(\bm{q}\right) = &\lVert{\bm{A}\bm{q}-\bm{b}}\rVert_2^2\\ 
 		= &\left(\bm{A}\bm{q}-\bm{b}\right)^{\mathsf{H}} \left(\bm{A}\bm{q}-\bm{b}\right)\\
 		=&
 		\bm{q}^{\mathsf{H}} \bm{A}^{\mathsf{H}} \bm{A} \bm{q} - \bm{q}^{\mathsf{H}} \bm{A}^{\mathsf{H}} \bm{b} -\bm{b}^{\mathsf{H}} \bm{A} \bm{q} + \bm{b}^{\mathsf{H}} \bm{b}.
 	\end{aligned}
 \end{align} 
 Upon using the first, the second and the fourth rows of TABLE $\mathrm{\ref{table}}$, we take the gradient of $f\left(\bm{q}\right)$ with respect to $\boldsymbol{q}^*$ to gield
\begin{equation}\label{2}
	\begin{aligned}
		\nabla_{\bm{q}^*} f\left(\bm{q}\right)
		\triangleq &\left(\dfrac{\partial f}{\partial \bm{q}^*}\right)^{\mathsf{T}}
		\overset{\eqref{cong}}{=} \left(\dfrac{\partial f}{\partial \bm{q}}\right)^{\mathsf{H}}\\
		= &\dfrac{1}{2} \bm{A}^{\mathsf{H}} \bm{A}\bm{q} + \dfrac{1}{2} \bm{A}^{\mathsf{H}} \bm{b} - \bm{A}^{\mathsf{H}} \bm{b}\\
		= &\dfrac{1}{2} \bm{A}^{\mathsf{H}}  \left(\bm{A}\bm{q}-\bm{b}\right).
	\end{aligned}
\end{equation}
Then $\forall \bm{p}, \bm{q} \in \mathbb{H}^n$, we obtain
\begin{align}
	\begin{aligned}
		&f\left(\bm{q}\right) - f\left(\bm{p}\right) - 4{\rm Re}\left\{\nabla_{\bm{p}^*}f\left(\bm{p}\right)^{\mathsf{H}}\left(\bm{q}-\bm{p}\right)\right\}
		\\=&
		\left(\bm{A}\bm{q}-\bm{b}\right)^{\mathsf{H}}\left(\bm{A}\bm{q}-\bm{b}\right) - \left(\bm{A}\bm{p}-\bm{b}\right)^{\mathsf{H}}\left(\bm{A}\bm{p}-\bm{b}\right) 
		\\-& 2{\rm Re}\left\{\left(\bm{A}^{\mathsf{H}} \left(\bm{A}\bm{p}-\bm{b}\right)\right)^{\mathsf{H}}\left(\bm{q}-\bm{p}\right)\right\}
%			\\=&
%			\bm{q}^{\mathsf{H}} \bm{A}^{\mathsf{H}} \bm{A} \bm{q} - \bm{q}^{\mathsf{H}} \bm{A}^{\mathsf{H}} \bm{b} - \bm{b}^{\mathsf{H}} \bm{A} \bm{q} - \bm{p}^{\mathsf{H}} \bm{A}^{\mathsf{H}} \bm{A} \bm{p} + \bm{p}^{\mathsf{H}} \bm{A}^{\mathsf{H}} \bm{b} + \bm{b}^{\mathsf{H}} \bm{A} \bm{p} \nonumber \\&
%			- \big(\bm{p}^{\mathsf{H}} \bm{A}^{\mathsf{H}} \bm{A} \bm{q} + \bm{q}^{\mathsf{H}} \bm{A}^{\mathsf{H}} \bm{A} \bm{p} - 2\bm{p}^{\mathsf{H}} \bm{A}^{\mathsf{H}} \bm{A} \bm{p} - \bm{b}^{\mathsf{H}} \bm{A} \bm{q} - \bm{q}^{\mathsf{H}} \bm{A}^{\mathsf{H}} \bm{b}\\ +& \bm{b}^{\mathsf{H}} \bm{A} \bm{p} + \bm{p}^{\mathsf{H}} \bm{A}^{\mathsf{H}} \bm{b} \big)\nonumber
		\\=&
		\bm{q}^{\mathsf{H}} \bm{A}^{\mathsf{H}} \bm{A} \bm{q} + \bm{p}^{\mathsf{H}} \bm{A}^{\mathsf{H}} \bm{A} \bm{p} -\bm{p}^{\mathsf{H}} \bm{A}^{\mathsf{H}} \bm{A} \bm{q} - \bm{q}^{\mathsf{H}} \bm{A}^{\mathsf{H}} \bm{A} \bm{p}
		\\=&
		\left(\bm{q}-\bm{p}\right)^{\mathsf{H}} \bm{A}^{\mathsf{H}} \bm{A} \left(\bm{q} - \bm{p}\right)
		\\=&
		\lVert{\bm{A} \left(\bm{q} - \bm{p}\right)}\rVert_2^2
		\\\geqslant & 0.
	\end{aligned}
\end{align}
%where the last step uses 
%%	$\forall \bm{A} \in \mathbb{H}^{m\times n}$,
%$\bm{A}^{\mathsf{H}} \bm{A}\succeq  \bm{O}$.
Therefore, from Theorem $\ref{yijie}$ we know that $f\left(\bm{q}\right)$ is convex.\\

(\textbf{Gradient monotonicity criterion}) 
%, applying Theorem $\ref{dandiao}$)
$\forall \bm{p}, \bm{q} \in \mathbb{H}^n$, we have
\begin{align}
	\begin{aligned}
		&{\rm Re}\left\{\left(\nabla_{\bm{p}^*}f\left(\bm{p}\right)-\nabla_{\bm{q}^*}f\left(\bm{q}\right)\right)^{\mathsf{H}}\left(\bm{p}-\bm{q}\right)\right\}
		\\=&\dfrac{1}{2} {\rm Re}\left\{\left(\bm{A}^{\mathsf{H}} \left(\bm{A}\bm{p}-\bm{b}\right) - \bm{A}^{\mathsf{H}} \left(\bm{A}\bm{q}-\bm{b}\right)\right)^{\mathsf{H}}\left(\bm{p}-\bm{q}\right)\right\}
		\\=&
		\dfrac{1}{2} {\rm Re}\left\{\left(\bm{p}-\bm{q}\right)^{\mathsf{H}} \bm{A}^{\mathsf{H}} \bm{A} \left(\bm{p} - \bm{q}\right)\right\}
		\\=&
		\dfrac{1}{2} \left(\bm{p}-\bm{q}\right)^{\mathsf{H}} \bm{A}^{\mathsf{H}} \bm{A} \left(\bm{p} - \bm{q}\right)
		\\=&
		\dfrac{1}{2} \lVert{\bm{A} \left(\bm{q} - \bm{p}\right)}\rVert_2^2
		\\\geqslant & 0.
	\end{aligned}
\end{align}
Therefore, from Theorem $\ref{dandiao}$ we know that $f\left(\bm{q}\right)$ is convex.\\

(\textbf{Second-order characterization criterion})
%	(??4??????)U
Using the third row of TABLE $\mathrm{\ref{table}}$, we get
\begin{align}
	\begin{aligned}\label{hqq**}
		\boldsymbol{H}_{\bm{q} \bm{q}^{*}} 
		\triangleq &
		\dfrac{\partial}{\partial \bm{q}}\left(\dfrac{\partial f}{\partial \bm{q}^{*}}\right)^{\mathsf{T}} 
		= \dfrac{\partial \nabla_{\bm{q}^*} f\left(\bm{q}\right)}{\partial \bm{q}}
		\\\overset{\eqref{2}}{=} &
		\dfrac{\partial \left(\dfrac{1}{2} \bm{A}^{\mathsf{H}} \left(\bm{A}\bm{q}-\bm{b}\right)\right)}{\partial \bm{q}}
		=
		\dfrac{1}{2}\bm{A}^{\mathsf{H}} \bm{A},
	\end{aligned}
\end{align}
and for any $\nu \in \left\{i,j,k \right\}$, 
\begin{align}\label{hvq*}
	\begin{aligned}
		\boldsymbol{H}_{\bm{q}^\nu \bm{q}^{*}} 
		\triangleq &
		\dfrac{\partial}{\partial \bm{q}^{\nu}}\left(\dfrac{\partial f}{\partial \bm{q}^{*}}\right)^{\mathsf{T}} 
		= \dfrac{\partial \nabla_{\bm{q}^*} f\left(\bm{q}\right)}{\partial \bm{q}^{\nu}}\\
		\overset{\eqref{2}}{=} &
		\dfrac{\partial \left(\dfrac{1}{2} \bm{A}^{\mathsf{H}} \left(\bm{A}\bm{q}-\bm{b}\right)\right)}{\partial \bm{q}^{\nu}}
		\overset{\eqref{pr}}{=}  
		\dfrac{1}{2} \bm{A}^{\mathsf{H}}\bm{A} \dfrac{\partial \bm{q}}{\partial \bm{q}^{\nu}}
		=\bm{O}.
	\end{aligned}
\end{align}
Then
$\forall \bm{x} \in \mathbb{H}^n, \bm{x} \ne \bm{0}$, it follows that
\begin{align}
	\begin{aligned}
		&\sum\limits_{\nu \in \left\{1,i,j,k \right\}} {\rm Re}\left\{\bm{x}^{\mathsf{H}} \boldsymbol{H}_{\bm{q}^\nu \bm{q}^{*}} \bm{x}^{\nu}\right\}\\
		=&
		\dfrac{1}{2} {\rm Re}\left\{\bm{x}^{\mathsf{H}} \bm{A}^{\mathsf{H}} \bm{A} \bm{x}\right\}
		=
		\dfrac{1}{2} \bm{x}^{\mathsf{H}} \bm{A}^{\mathsf{H}} \bm{A} \bm{x}
		=
		\dfrac{1}{2} \lVert{\bm{A} \bm{x}}\rVert_2^2
		\geqslant 0
		.
	\end{aligned}
\end{align}
Therefore, by Corollary $\ref{ertuilun}$, we know that $f\left(\bm{q}\right)$ is convex.
\end{IEEEproof}
%\begin{example}
%	If the quaternion function 
%	$f\left(\bm{q}\right)=\lVert{\bm{A}\bm{q}-\bm{b}}\rVert_2^2$, $\forall \bm{q} \in \mathbb{H}^n$, $\bm{b}\in \mathbb{H}^m$, the matrix $\bm{A}\in \mathbb{H}^{m\times n}$ is nonsingular, then $f\left(\bm{q}\right)$ is strictly convex.
%	%????$f\left(\bm{q}\right)=\lVert{\bm{A}\bm{q}-\bm{b}}\rVert_2^2$, $ \forall \bm{q} \in \mathbb{H}^n$, ????? $\bm{A}\in \mathbb{H}^{m\times n}$, $\bm{b}\in \mathbb{H}^m$, ?$f\left(\bm{q}\right)$??????.
%\end{example}
%\begin{IEEEproof}
%The proof is similar to that of Example $\ref{aq-b}$.
%\end{IEEEproof}

\section{Strongly Convex Quaternion Function: Definition and Discriminant Theorems}
We shall now discuss the discriminant criteria for strongly convex quaternion functions, building upon the theorems for convexity. These criteria will be useful in designing optimization algorithms. %that take advantage of the additional structure provided by strong convexity.
\begin{definition}[Strongly convex function]\label{qtdj}
	%(The equivalent definition of strongly convex function)
	The quaternion function $f\left(\bm{q}\right):\mathcal{C} \subset \mathbb{H}^n \to \mathbb{R}$ is called strongly convex, if $ \exists \sigma>0 $, $ \forall \bm{p}$, $\bm{q} \in \mathcal{C} $, $\forall \theta \in(0,1) $, 
	\begin{equation}\label{qtt1}
		f\big(\theta \bm{p}+(1-\theta) \bm{q}\big) \leqslant \theta f(\bm{p})+(1-\theta) f(\bm{q})-\dfrac{\sigma}{2} \theta(1-\theta)\|\bm{p}-\bm{q}\|^2_2,
	\end{equation}
	where $\sigma$ is the strongly convex parameter. For convenience, $f\left(\bm{q}\right)$ is also called $\sigma$-strongly convex.
\end{definition}

Based on the definition of strongly convex functions, we obtain the following equivalence theorem.
%\begin{theorem}
%	The quaternion function $f\left(\bm{q}\right):\mathcal{C} \subset \mathbb{H}^n \to \mathbb{R}$ has a unique minimum value; if $f\left(\bm{q}\right)$ is strongly convex and has a minimum value.
%\end{theorem}
%\begin{IEEEproof} 
%The proof is the same as that in the real field \cite{boyd_vandenberghe_2004,nesterov2018lectures}.
%\end{IEEEproof}
\begin{theorem}\label{qt}
	The quaternion function $f\left(\bm{q}\right):\mathcal{C} \subset \mathbb{H}^n \to \mathbb{R}$ is $\sigma$-strongly convex, iff $\exists \sigma>0$, s.t. the function
	\begin{equation}
		g\left(\bm{q}\right) \triangleq f\left(\bm{q}\right) - \dfrac{\sigma}{2} \lVert{\bm{q}}\rVert^2_2 
	\end{equation}
	is convex.
\end{theorem}
\begin{IEEEproof} 
This is straightforward to prove, by applying Definition $\ref{tf}$ and Definition $\ref{qtdj}$.
\end{IEEEproof}

Similar to convex quaternion functions, strongly convex quaternion functions also have first-order characterization, gradient monotonicity, and second-order characterization.
\begin{theorem}
	%(Quadratic lower bound)\textcolor{red}{????????}
	[First-order characterization]
	Consider a convex set $\mathcal{C}\subset \mathbb{H}^n$ and a real-differentiable quaternion function $f(\bm{q}):\mathcal{C} \to \mathbb{R}$. Then, $f(\bm{q})$ is $\sigma$-strongly convex iff $\forall \bm{p}, \bm{q} \in \mathcal{C}$,
	%(????)???$f\left(\bm{q}\right):\mathcal{C} \subset \mathbb{H}^n \to \mathbb{R}$?????$m$-????, ????????: 
	\begin{equation}\label{ecxjgs}
		f\left(\bm{q}\right) \geqslant
		f\left(\bm{p}\right) + 4{\rm Re}\left\{\nabla_{\bm{p}^*} f\left(\bm{p}\right)^{\mathsf{H}} \left(\bm{q}-\bm{p}\right)\right\} + \dfrac{\sigma}{2}\lVert{\bm{q}-\bm{p}}\rVert_2^2,
	\end{equation}
	where $\nabla_{\bm{p}^*}f\left(\bm{p}\right)$ is defined in $(\ref{getd})$.
\end{theorem}
\begin{IEEEproof} 
	From Theorem $\ref{qt}$, $f(\bm{p})$ is strongly convex iff $g\left(\bm{p}\right) = f\left(\bm{p}\right) - \frac{1}{2}\sigma \lVert{\bm{p}}\rVert^2_2$ is convex. Then, upon applying Theorem $\ref{yijie}$, $\forall \bm{p}, \bm{q} \in \mathcal{C}$,
	\begin{equation}
		g\left(\bm{q}\right) \geqslant g\left(\bm{p}\right)+4{\rm Re}\left\{\nabla_{\bm{p}^*}g\left(\bm{p}\right)^{\mathsf{H}}\left(\bm{q}-\bm{p}\right)\right\}.
	\end{equation}
	Using the fourth row of TABLE $\mathrm{\ref{table}}$, $\nabla_{\bm{p}^*} g\left(\bm{p}\right) = \nabla_{\bm{p}^*}f\left(\bm{p}\right) - \frac{1}{4}\sigma\bm{p}$. Then, $\forall \bm{p}, \bm{q} \in \mathcal{C}$,
	\begin{align}\label{qt3}
		&f\left(\bm{q}\right) - \dfrac{\sigma}{2} \lVert{\bm{q}}\rVert^2_2 \\\geqslant& f\left(\bm{p}\right) - \dfrac{\sigma}{2} \lVert{\bm{p}}\rVert^2_2 + 4{\rm Re}\left\{\big(\nabla_{\bm{p}^*}f\left(\bm{p}\right)-\dfrac{\sigma}{4}\bm{p}\big)^{\mathsf{H}}\left(\bm{q}-\bm{p}\right)\right\}.\nonumber
	\end{align}
	Since
	\begin{align}\label{1jpq}
		&\dfrac{\sigma}{2} \lVert{\bm{q}}\rVert^2_2 - \dfrac{\sigma}{2} \lVert{\bm{p}}\rVert^2_2 + 4{\rm Re}\left\{\big(\nabla_{\bm{p}^*}f\left(\bm{p}\right)-\dfrac{\sigma}{4}\bm{p}\big)^{\mathsf{H}}\left(\bm{q}-\bm{p}\right)\right\} \nonumber
		\\=&
		4{\rm Re}\left\{\nabla_{\bm{p}^*} f\left(\bm{p}\right)^{\mathsf{H}} \left(\bm{q}-\bm{p}\right)\right\} - \sigma{\rm Re}\left\{\bm{p}^{\mathsf{H}}\bm{q}\right\} + \sigma \lVert{\bm{p}}\rVert^2_2 \nonumber\\+ &\dfrac{\sigma}{2} \lVert{\bm{q}}\rVert^2_2 - \dfrac{\sigma}{2} \lVert{\bm{p}}\rVert^2_2\nonumber
		\\\overset{\eqref{sq}}{=}&
		4{\rm Re}\left\{\nabla_{\bm{p}^*} f\left(\bm{p}\right)^{\mathsf{H}} \left(\bm{q}-\bm{p}\right)\right\} + \dfrac{\sigma}{2}\lVert{\bm{q}-\bm{p}}\rVert_2^2.
	\end{align}
	Upon substituting $(\ref{1jpq})$ into $(\ref{qt3})$, the proof follows.
\end{IEEEproof}
\begin{theorem}[Gradient monotonicity]\label{mqt}
	Consider a convex set $\mathcal{C}\subset \mathbb{H}^n$ and a real-differentiable quaternion function $f(\bm{q}):\mathcal{C} \to \mathbb{R}$. Then, 
	$f(\bm{q})$ is $\sigma$-strongly convex iff $\forall \bm{p}, \bm{q} \in \mathcal{C}$,
	%?????$f(\bm{q}):\mathcal{C} \to \mathbb{R}$??????, ??$\mathcal{C}\subset \mathbb{H}^n$???, ?$ f $?$m$-????????
	\begin{equation}\label{qtu}
		{\rm Re}\left\{\big(\nabla_{\bm{p}^*} f\left(\bm{p}\right)-\nabla_{\bm{q}^*} f\left(\bm{q}\right)\big)^{\mathsf{H}}\left(\bm{p}-\bm{q}\right)\right\} \geqslant \dfrac{\sigma}{4} \lVert{\bm{p}-\bm{q}}\rVert_2^2, 
	\end{equation}
	where $\nabla_{\bm{p}^*}f\left(\bm{p}\right)$ is defined in $(\ref{getd})$.
\end{theorem}
\begin{IEEEproof} 
From the Theorem $\ref{qt}$, $f(\bm{q})$ is strongly convex iff $g\left(\bm{q}\right) = f\left(\bm{q}\right) - \frac{1}{2}\sigma \lVert{\bm{q}}\rVert^2_2$ is convex. Then, after applying Theorem $\ref{dandiao}$, we have
\begin{equation}
	{\rm Re}\left\{\left(\nabla_{\bm{p}^*} g\left(\bm{p}\right)-\nabla_{\bm{q}^*} g\left(\bm{q}\right)\right)^{\mathsf{H}}\left(\bm{p}-\bm{q}\right)\right\} \geqslant 0, ~ \forall \bm{p}, \bm{q} \in \mathcal{C}.
\end{equation}
Using the fourth row of TABLE $\mathrm{\ref{table}}$, we have
$\nabla_{\bm{q}^*} g\left(\bm{q}\right) = \nabla_{\bm{q}^*}f\left(\bm{q}\right) - \frac{1}{4}\sigma\bm{q}$, then $\forall \bm{p}, \bm{q} \in \mathcal{C}$,
\begin{equation}\label{qt2}
	{\rm Re}\left\{\big(\nabla_{\bm{p}^*} f\left(\bm{p}\right) - \dfrac{\sigma}{4}\bm{p} - \nabla_{\bm{q}^*}f\left(\bm{q}\right) + \dfrac{\sigma}{4}\bm{q}\big)^{\mathsf{H}}\left(\bm{p}-\bm{q}\right)\right\} \geqslant 0.
\end{equation}
Upon rearranging the terms in $(\ref{qt2})$, we obtain $(\ref{qtu})$. 
\end{IEEEproof}
\begin{theorem}
	[Second-order characterization]\label{4.5}
	Consider a convex set $\mathcal{C}\subset \mathbb{H}^n$ and a second-order continuous real-differentiable quaternion function $f(\bm{q}):\mathcal{C} \to \mathbb{R}$. Then, $f(\bm{q})$ is $\sigma$-strongly convex, iff 
	%(????)?????$f(\bm{q}):\mathcal{C} \to \mathbb{R}$??????????, ??$\mathcal{C}\subset \mathbb{H}^n$???, ? $ f $ ????????
	\begin{equation}
		\boldsymbol{H}_{\mathcal{H}\mathcal{H}^*} \succeq \dfrac{\sigma}{4}\bm{I}_{4n},
	\end{equation}
	where $\boldsymbol{H}_{\mathcal{H}\mathcal{H}^*}$ is defined in $(\ref{HHH})$.
\end{theorem}
\begin{IEEEproof} 
	Define $g\left(\bm{q}\right) \triangleq f\left(\bm{q}\right) -\frac{1}{2}\sigma \lVert{\bm{q}}\rVert^2_2$, $h(\bm{q}) \triangleq 2 \lVert{\bm{q}}\rVert^2_2 = 2\bm{q}^{\mathsf{H}}\bm{q} = 2\big(\bm{q}^{\mathsf{H}}\bm{q}\big)^{\mu} \overset{\eqref{xzp}}{=} 2\bm{q}^{\mu\mathsf{H}}\bm{q}^{\mu}$, $\mu \in \left\{1,i,j,k \right\}$. Upon applying the fourth row of TABLE $\mathrm{\ref{table}}$, we have
	\begin{align}\label{93}
		\left(\dfrac{\partial h}{\partial \bm{q}^{\mu *}}\right)^{\mathsf{T}} 
		= \bm{q}^{\mu}, \quad \mu \in \left\{1,i,j,k \right\}.
	\end{align}
	Then, $\forall \mu \in \left\{1,i,j,k \right\}$,
	\begin{align}
%		\boldsymbol{H}_{\bm{q}^\mu \bm{q}^{\mu*}} 
%			\triangleq 
		\dfrac{\partial}{\partial \bm{q}^{\mu}}\left(\dfrac{\partial h}{\partial \bm{q}^{\mu*}}\right)^{\mathsf{T}} 
			\overset{\eqref{93}}{=} 
		\dfrac{\partial \bm{q}^{\mu}}{\partial \bm{q}^{\mu}}
			= \bm{I}_{n},
	\end{align}
	and $\forall \mu, \nu \in \left\{1,i,j,k \right\}$, $\mu \ne \nu$,
	\begin{align}
%		\boldsymbol{H}_{\bm{q}^\nu \bm{q}^{\mu*}} 
%		\triangleq 
		\dfrac{\partial}{\partial \bm{q}^{\nu}}\left(\dfrac{\partial h}{\partial \bm{q}^{\mu*}}\right)^{\mathsf{T}} 
		\overset{\eqref{93}}{=} 
		\dfrac{\partial \bm{q}^{\mu}}{\partial \bm{q}^{\nu}}
		= \bm{O}.
	\end{align}
	By $(\ref{HHH})$, the augmented quaternion Hessian matrix of $h$ is $\bm{I}_{4n}$. Therefore, the augmented quaternion Hessian matrix of $g$ is $\boldsymbol{H}_{\mathcal{H}\mathcal{H}^*} -\frac{1}{4}\sigma \bm{I}_{4n}$.
	From Theorem $\ref{qt}$, $f(\bm{q})$ is strongly convex iff $g\left(\bm{q}\right)$ is convex. 
	Then upon applying Theorem $\ref{er}$, $g\left(\bm{q}\right) = f\left(\bm{q}\right) -\frac{1}{2}\sigma \lVert{\bm{q}}\rVert^2_2$ is convex iff  $\boldsymbol{H}_{\mathcal{H}\mathcal{H}^*} - \frac{1}{4}\sigma \bm{I}_{4n} \succeq \bm{O}$.
\end{IEEEproof}
\begin{corollary}
	Consider a convex set $\mathcal{C}\subset \mathbb{H}^n$ and a second-order continuous real-differentiable quaternion function $f(\bm{q}):\mathcal{C} \to \mathbb{R}$. Then, the following three propositions are equivalent:
	%	?????$f(\bm{q}):\mathcal{C} \to \mathbb{R}$??????????, ??$\mathcal{C}\subset \mathbb{H}^n$???, ?????????
	
	$\mathit{(a)}$~ $ f(\bm{q}) $ is $\sigma$-strongly convex;
	
	$\mathit{(b)}$~ $\boldsymbol{H}_{\mathcal{H}\mathcal{H}^*}\succeq  \frac{1}{4}\sigma \bm{I}_{4n}$;
	
	$\mathit{(c)}$~ $\sum\limits_{\nu \in \left\{1,i,j,k \right\}} {\rm Re}\left\{\bm{x}^{\mathsf{H}} \boldsymbol{H}_{\bm{q}^\nu \bm{q}^{*}} \bm{x}^{\nu}\right\} - \frac{1}{4}\sigma \lVert{\bm{x}}\rVert^2_2
	\geqslant 0
	$,\quad $\forall \bm{x} \in \mathbb{H}^n$, $\bm{x} \ne \bm{0} .$
\end{corollary}
\begin{IEEEproof} 
	According to Theorem
	$\ref{4.5}$, $(a)$ is equivalent to $(b)$, so we only need to prove that $(b)$ is equivalent to $(c)$. 
%	Since from Corollary $\ref{ertuilun}$,
	By Corollary $\ref{hermite}$,  $\boldsymbol{H}_{\mathcal{H}\mathcal{H}^*}$ is a Hermite matrix, so $\boldsymbol{H}_{\mathcal{H}\mathcal{H}^*} - \frac{1}{4}\sigma \bm{I}_{4n}$ is also Hermite matrix. Then, $\forall \bm{x}_{\mathcal{H}} \in \mathcal{H}$, $\bm{x}_{\mathcal{H}} \ne \bm{0}$, we have
	%??$(\ref{Hermite})$,$\boldsymbol{H}_{\mathcal{H}\mathcal{H}^*}$?Hermite????, ?
	\begin{align}
	\begin{aligned}
			&\quad~\bm{x}_{\mathcal{H}} ^{\mathsf{H}} \left(\boldsymbol{H}_{\mathcal{H}\mathcal{H}^*}- \dfrac{\sigma}{4} \bm{I}_{4n}\right) \bm{x}_{\mathcal{H}}\\
		&= \bm{x}_{\mathcal{H}} ^{\mathsf{H}} \boldsymbol{H}_{\mathcal{H}\mathcal{H}^*} \bm{x}_{\mathcal{H}} -\dfrac{\sigma}{4} \bm{x}_{\mathcal{H}} ^{\mathsf{H}} \bm{x}_{\mathcal{H}}
		\\&\overset{\eqref{fanshu}}{=}
		%		\sum_{\nu,\mu \in \left\{1,i,j,k \right\}} \bm{x}^{\mu \mathsf{H}} \dfrac{\partial}{\partial \bm{q}^\nu}\left(\dfrac{\partial f}{\partial \bm{q}^{\mu *}}\right)^{\mathsf{T}} \bm{x}^{\nu}
		%		=
		\sum_{\mu,\nu \in \left\{1,i,j,k \right\}} \bm{x}^{\mu \mathsf{H}} \boldsymbol{H}_{\bm{q}^\nu \bm{q}^{\mu*}} \bm{x}^{\nu}- \sigma \bm{x}^{\mathsf{H}} \bm{x}
		\\&\overset{\eqref{r,h}}{=}
		%		4\sum_{\mu \in \left\{1,i,j,k \right\}} {\rm Re}\left(\bm{x}^{\mu \mathsf{H}} \boldsymbol{H}_{\bm{q} \bm{q}^{\mu*}} \bm{x}\right)
		%		=
		4\sum_{\nu \in \left\{1,i,j,k \right\}} {\rm Re}\left\{\bm{x}^{\mathsf{H}} \boldsymbol{H}_{\bm{q}^\nu \bm{q}^{*}} \bm{x}^{\nu}\right\}-\sigma \lVert{\bm{x}}\rVert^2_2.
	\end{aligned}
	\end{align}
	%	\left(??
	%	\begin{equation}
		%		\bm{x}_{\mathcal{H}} ^{\mathsf{H}} \boldsymbol{H}_{\mathcal{H}\mathcal{H}^*} \bm{x}_{\mathcal{H}} 
		%		= 
		%		16\bm{x}^{\mathsf{H}} \boldsymbol{H}_{\bm{q}\bm{q}^*} \bm{x}
		%	\end{equation}\right)
	Therefore, $\forall \bm{x} \in \mathbb{H}^n$, $\bm{x} \ne \bm{0}$,
	\begin{align}
		\begin{aligned}
			\sum_{\nu \in \left\{1,i,j,k \right\}} &{\rm Re}\left\{\bm{x}^{\mathsf{H}} \boldsymbol{H}_{\bm{q}^\nu \bm{q}^{*}} \bm{x}^{\nu}\right\}-\dfrac{\sigma}{4} \lVert{\bm{x}}\rVert^2_2\geqslant 0, \\ \quad &\Leftrightarrow \quad \boldsymbol{H}_{\mathcal{H}\mathcal{H}^*} - \dfrac{\sigma}{4} \bm{I}_{4n} \succeq \bm{O}.
		\end{aligned}
	\end{align}
This completes the proof.
\end{IEEEproof}

Upon applying Lemma $\ref{zhuzishi}$, we can obtain a necessary condition for $\sigma$-strongly convex quaternion functions.
\begin{theorem}
	Consider a convex set $\mathcal{C}\subset \mathbb{H}^n$ and a second-order continuous real-differentiable quaternion function $f(\bm{q}):\mathcal{C} \to \mathbb{R}$. If $f(\bm{q})$ is $\sigma$-strongly convex, then
	%the quaternion Hessian matrix $\boldsymbol{H}_{\bm{q}\bm{q}^*}$ is positive semi-definite.
	\begin{equation}
		\boldsymbol{H}_{\bm{q}\bm{q}^*} \succeq \dfrac{\sigma}{4} \bm{I}_{n},
	\end{equation}
	%	In addition, $f(\bm{q})$ is strictly convex if 
	%	%$ f $ ??????????
	%	\begin{equation}
		%		\boldsymbol{H}_{\mathcal{H}\mathcal{H}^*} \succ \bm{O},
		%	\end{equation}
	where $\boldsymbol{H}_{\bm{q}\bm{q}^*}$ is the quaternion Hessian matrix, defined in $(\ref{hqq*})$.
	%?????$f(\bm{q}):\mathcal{C} \to \mathbb{R}$??????????, ??$\mathcal{C}\subset \mathbb{H}^n$???, ? $ f $ ????, ?$\boldsymbol{H}_{\bm{q}\bm{q}^*}$???????? $ f $ ??????, ?$\boldsymbol{H}_{\bm{q}\bm{q}^*}$?????.
\end{theorem}
\begin{IEEEproof} 
	Note that $f(\bm{q})$ is $\sigma$-strongly convex, and upon applying Theorem $\ref{4.5}$, we have $\boldsymbol{H}_{\mathcal{H}\mathcal{H}^*} \succeq \frac{1}{4}\sigma \bm{I}_{4n}$. By $(\ref{HHH})$ and Lemma $\ref{zhuzishi}$, we finally obtain 
	$\boldsymbol{H}_{\bm{q}\bm{q}^*} \succeq \frac{1}{4}\sigma \bm{I}_{n}$.
	%	???$\ref{er}$???, $ f $ ???????$\boldsymbol{H}_{\mathcal{H}\mathcal{H}^*}\succeq  \bm{O}$.????$\ref{zhengding}$?, $\boldsymbol{H}_{\bm{q}\bm{q}^*}$??????.\quad $ f $ ???????, ????.
\end{IEEEproof}

\section{Convex Quaternion Optimization Problems and Their Applications in Signal Processing}
We now proceed to introduce the convex quaternion problem and its fundamental theorem. This is followed by several applications of convex quaternion optimization in communications, highlighting its practical significance.
\subsection{Convex Quaternion Optimization Problems}
Similar to convex real and complex optimization problems, 
convex quaternion optimization problems generally have a structure which consist of the minimization of a convex quaternion function subject to (shortened to s.t.) quaternion affine equality constraints and inequality constraints defined by convex quaternion functions, as follows
%????????????????????:
\begin{align}\label{qtyh}
	\begin{aligned}
		\min \limits_{\bm{q} \in \mathbb{H}^n}& f_{0}(\bm{q}) \\
		\text { s.t. } & \bm{A q}=\bm{b},\\ 
		& f_{i}(\bm{q}) \leqslant 0, ~i=1,\, \ldots,\, P
	\end{aligned}
\end{align}
where $ f_{i}: \mathbb{H}^n \to \mathbb{R}$, $i=0,\, 1,\, \ldots, \,P $ is convex, $ \bm{A} \in \mathbb{H}^{m \times n}$, $\bm{b} \in \mathbb{H}^{m} $. The problem field is 
$\mathcal{F} \triangleq \bigcap_{i=0}^{P} \bm{dom} f_{i}$, 
and feasible set is 
$\mathcal{C} \triangleq$ $\left\{\bm{q} \in \mathcal{F} \mid f_{i}(\bm{q}) \leqslant 0,~ i=1,\, \ldots,\, P,~ \bm{A} \bm{q}=\bm{b}\right\}$. 
From Definition $\ref{tf}$, Example $\ref{aq=b}$ and Example $\ref{fq}$, the sets $\mathcal{D} \triangleq \left\{\bm{q} \in \mathbb{H}^{n} \mid  \bm{A} \bm{q} = \bm{b}\right\}$, $\mathcal{E}_i \triangleq \left\{\bm{q} \in \mathbb{H}^{n} \mid  f_i(\bm{q}) \leqslant 0\right\}$, $i=1,\, \ldots,\, P$ and $\bm{dom} f_{i}$, $i=0,\,1,\, \ldots,\, P$
are convex. Therefore, the set $\mathcal{C} = \mathcal{D} \bigcap \left(\bigcap_{i=1}^{P} \mathcal{E}_i\right) \bigcap \mathcal{F}$ is also convex. 

%For example, the quaternion least squares problem:
%\begin{equation}
%	\min \limits_{\bm{q} \in \mathbb{H}^n} \|\bm{A q}-\bm{b}\|_{2}^{2}
%\end{equation}
%is an unconstrained convex quaternion optimization problem, because  $f\left(\bm{q}\right)=\lVert{\bm{A}\bm{q}-\bm{b}}\rVert_2^2$ is convex, from Example 
%$\ref{aq-b}$.

When studying the convexity of quaternion functions, we utilized the augmented quaternion vectors and augmented real vectors. Similarly, when studying the properties of quaternion convex optimization problems, we also need to utilize the augmented quaternion and the augmented real convex optimization settings.

%????????$(\ref{qtyh})$???????????$(\ref{htyh})$?????????$(\ref{rtyh})$????.\\
%??????????:
The convex augmented quaternion optimization problem of $(\ref{qtyh})$ is given by \cite{flamant2021general}
\begin{align}\label{htyh}
	\begin{aligned}
		\min \limits_{\bm{q}_\mathcal{H} \in \mathcal{H}} & f_{0}(\bm{q}_\mathcal{H}) \\
		\text { s.t. } & \bm{A}_\mathcal{H} \bm{q}_\mathcal{H}=\bm{b}_\mathcal{H},\\ & f_{i}(\bm{q}_\mathcal{H}) \leqslant 0,~ i=1,\, \ldots,\, P
	\end{aligned}
\end{align}
where $ f_{i}: \mathcal{H} \to \mathbb{R}$, $i=0,\, 1,\, \ldots,\, P $ is convex, $\bm{A}_\mathcal{H} \triangleq {\rm diag}\left(\bm{A},\,\bm{A}^i,\,\bm{A}^j,\,\bm{A}^k\right) \in \mathbb{H}^{4m \times 4n}$, $\bm{A} \in \mathbb{H}^{m \times n}$, and $ \bm{b}_\mathcal{H} \triangleq \left(\bm{b}^{\mathsf{T}},\bm{b}^{i\mathsf{T}},\bm{b}^{j\mathsf{T}},\bm{b}^{k\mathsf{T}}\right)^{\mathsf{T}} \in \mathbb{H}^{4m} $.

The convex augmented real convex optimization problem of $(\ref{qtyh})$ is given by \cite{flamant2021general}
%????????:
\begin{align}\label{rtyh}
	\begin{aligned}
		\min \limits_{\bm{q}_\mathcal{R} \in \mathcal{R}} & f_{0}(\bm{q}_\mathcal{R}) \\
		\text { s.t. } & \bm{A}_\mathcal{R} \bm{q}_\mathcal{R}=\bm{b}_\mathcal{R},\\ & f_{i}(\bm{q}_\mathcal{R}) \leqslant 0,~ i=1,\, \ldots,\, P
	\end{aligned}
\end{align}
where $ f_{i}: \mathcal{R} \to \mathbb{R}$, $i=0,\, 1,\, \ldots,\, P $ is convex, together with $\bm{A}_\mathcal{R} \triangleq \frac{1}{4} \bm{J}_{m}^{\mathsf{H}} \bm{A}_\mathcal{H} \bm{J}_{n} \in \mathbb{R}^{4m \times 4n}$, $\bm{A}_\mathcal{H} \in \mathbb{H}^{4m \times 4n}$, and $\bm{b}_\mathcal{R} \triangleq \frac{1}{4} \bm{J}_m^{\mathsf{H}} \bm{b}_{\mathcal{H}} \in \mathbb{R}^{4m} $.
\begin{lemma}[\cite{flamant2021general}]\label{qhrdj}
	The convex quaternion optimization problem in $(\ref{qtyh})$, the convex augmented quaternion optimization problem in $(\ref{htyh})$, and the convex augmented real optimization problem in $(\ref{rtyh})$ are equivalent.
\end{lemma}
%\begin{IEEEproof} 
%From $(\ref{iff})$, we have
%\begin{equation}
%	f_{i}(\bm{q}_\mathcal{H}) = f_{i}(\bm{q}_\mathcal{R}) = f_{i}(\bm{q}),~ i=0,\,1,\, \ldots,\, P.
%\end{equation}
%Applying $(\ref{r,h})$, we obtain $\bm{A}_\mathcal{R} \bm{q}_\mathcal{R}=\frac{1}{16} \bm{J}_{m}^{\mathsf{H}} \bm{A}_\mathcal{H} \bm{J}_{n} \bm{J}_{n}^{\mathsf{H}} \bm{q}_\mathcal{H}=\frac{1}{4} \bm{J}_{m}^{\mathsf{H}} \bm{A}_\mathcal{H} \bm{q}_\mathcal{H}$, and $\bm{b}_\mathcal{R} = \frac{1}{4} \bm{J}_{m}^{\mathsf{H}} \bm{b}_\mathcal{H}$, then 
%\begin{equation}
%\bm{A}_\mathcal{R} \bm{q}_\mathcal{R}=\bm{b}_\mathcal{R} ~ \Leftrightarrow ~ 
%	\bm{A}_\mathcal{H} \bm{q}_\mathcal{H}=\bm{b}_\mathcal{H} ~ \Leftrightarrow ~ \bm{A} \bm{q}=\bm{b}.
%\end{equation}
%Therefore, we can conclude the proof.
%%??????????$(\ref{qtyh})$???????????$(\ref{htyh})$?????????$(\ref{rtyh})$??.
%\end{IEEEproof}
\begin{theorem}\label{qjzy}
	For the convex quaternion optimization problem in $(\ref{qtyh})$, any local optimal solution is also the global optimal solution.
	%??????????$(\ref{qtyh})$, ??????????????.
\end{theorem}
\begin{IEEEproof}
We already know \cite{2017Convexx} that for the real convex optimization problem in $(\ref{rtyh})$, any local optimal solution, for example $\bar{\bm{q}}_{\mathcal{R}}$, is also the global optimal solution. Then, from Lemma $\ref{qhrdj}$, the local optimal solution $\bar{\bm{q}}_{\mathcal{H}} \triangleq \left(\bar{\bm{q}}^{\mathsf{T}},\bar{\bm{q}}^{i\mathsf{T}},\bar{\bm{q}}^{j\mathsf{T}},\bar{\bm{q}}^{k\mathsf{T}}\right)^{\mathsf{T}} = \bm{J}_n\bar{\bm{q}}_{\mathcal{R}} $ is also global, in the augmented convex quaternion optimization problem in $(\ref{htyh})$. Therefore, the local optimal solution, $\bar{\bm{q}}$, is also global, in the convex quaternion optimization problem in $(\ref{qtyh})$.
\end{IEEEproof}
\subsection{Applications of Convex Quaternion Optimization in Signal Processing}
\begin{application}
	[Quaternion linear mean-square error filter]
%	A linear filter approximates the desired sequence $ d(n) $ through a linear combination of a window of input samples $ x(n) $ such that the estimate of the desired sequence is
%	
%
%	
%	where the input vector at time $ n $ is written as $ \bm{x}(n)=[x(n) x(n-1) \cdots x(n-N+1)]^{T} $ and the filter weights as $ \bm{w}=\left[w_{0} w_{1} \cdots w_{N-1}\right]^{T} $. 
%	
	The quaternion minimum mean-square error (MSE) filter can be specified as
	\begin{align}\label{AAAAAAA}
		\min\limits_{\bm{w} \in \mathbb{H}^n} \mathcal{J}(\bm{w}) 
		\triangleq
		E\left\{|e(n)|^{2}\right\}=E\left\{|d(n)-y(n)|^{2}\right\},
	\end{align}
	where $y(n)=\bm{w}^{\mathsf{H}} \bm{x}(n)$, $ \bm{x}(n) \in \mathbb{H}^n$ is the input vector, $ \bm{w} \in \mathbb{H}^n$ is the filter weight vector, and $ d(n)\in \mathbb{H}$ is the desired sequence. By the definition of the modulus, we have
	%From Example $\ref{aq-b}$, and 
	\begin{equation}\label{EEEEE}
		\begin{aligned}
			\mathcal{J}(\bm{w}) 
%			=
%			&E\left\{|d(n)-y(n)|^{2}\right\}\\
			=
			&E\left\{|d(n)-\bm{w}^{\mathsf{H}} \bm{x}(n)|^{2}\right\}\\
			=
%			&E\left\{|\bm{x}^{\mathsf{H}}(n)\bm{w} -d^*(n)|^{2}\right\},
			&E\left\{\left(d(n)-\bm{w}^{\mathsf{H}}\bm{x}(n)\right) \left(d(n)-\bm{w}^{\mathsf{H}}\bm{x}(n)\right)^*\right\}\\
			=
			&\bm{w}^{\mathsf{H}}E\{\bm{x}(n)\bm{x}^{\mathsf{H}}(n)\}\bm{w} - E\{d(n)\bm{x}^{\mathsf{H}}(n)\}\bm{w}\\ -& \bm{w}^{\mathsf{H}} E\{\bm{x}(n) d^*(n)\} + E\{d(n)d^*(n)\}\\
			=
			&\bm{w}^{\mathsf{H}}\bm{R}\bm{w} - \bm{p}^{\mathsf{H}}\bm{w} - \bm{w}^{\mathsf{H}} \bm{p} + \sigma_d^2
			,
		\end{aligned}
	\end{equation}
	where $\bm{R}=E\{\bm{x}(n)\bm{x}^{\mathsf{H}}(n)\}$ denotes the quaternion-valued input correlation matrix, $\bm{p}=E\{\bm{x}(n) d^*(n)\}$ is the crosscorrelation vector between the desired response and the input signal, $\sigma_d^2=E\{d(n)d^*(n)\}$ is the power of the desired response. By $(\ref{62})$ in Example $\ref{aq-b}$, we know that $\mathcal{J}(\bm{w})$ is convex. 
	Similarly to $(\ref{2})$, 
	%Using the first, second and fourth rows of TABLE $\mathrm{\ref{table}}$, 
	we take the gradient of $\mathcal{J}\left(\bm{w}\right)$ with respect to $\boldsymbol{w}^*$, and set the result to $\bm{0}$ to obtain
	\begin{equation}\label{JW}
		\begin{aligned}
			\nabla_{\boldsymbol{w}^*}\mathcal{J}\left(\bm{w}\right)
			=\left(\dfrac{\partial \mathcal{J}}{\partial \bm{w}}\right)^{\mathsf{H}}=
			\dfrac{1}{2} \bm{R} \bm{w} - \dfrac{1}{2}  \bm{p}
			= \bm{0}.
%			=&
%			E\left\{\dfrac{1}{2} \bm{x}(n) \bm{x}(n)^{\mathsf{H}} \bm{w} - \bm{x}(n) d(n)^* + \dfrac{1}{2} \bm{x}(n) d(n)^*\right\}\\
%			=
%			\dfrac{1}{2} E\left\{\bm{x}(n) \left(\bm{x}^{\mathsf{H}}(n) \bm{w} -  d(n)^*\right)\right\}.
		\end{aligned}
	\end{equation}
%	Then, the quaternion Hessian matrix of $\mathcal{J}\left(\bm{w}\right)$ is
%	\begin{equation}
%		\boldsymbol{H}_{\bm{q}^\mu \bm{q}^{\mu*}} 
%			\triangleq 
%		\dfrac{\partial}{\partial \bm{q}^{\mu}}\left(\dfrac{\partial \mathcal{J}}{\partial \bm{q}^{\mu*}}\right)^{\mathsf{T}} 
%	\end{equation}
%	Setting $\nabla_{\boldsymbol{w}^*}\mathcal{J}\left(\bm{w}\right)
%	= \bm{0}$, we have
%	\begin{align}
%		\begin{aligned}
%			\bm{R} \bm{w} = \bm{p}.
%%			&\quad~\nabla_{\boldsymbol{w}^*}\mathcal{J}\left(\bm{w}\right)
%%			= \bm{0} \\\quad
%%			\Rightarrow \quad 
%			%E\left\{\bm{x}(n) \bm{x}^{\mathsf{H}}(n)\right\} \bm{w}=E\left\{\bm{x}(n) d^{*}(n)\right\}.
%		\end{aligned}
%	\end{align}
	Using $(\ref{JW})$ and Theorem $\ref{qjzy}$, we arrive at gives the closed-form optimal solution
	% is obtained as
	\begin{equation}
		\bar{\bm{w}}=\bm{R}^{-1} \bm{p}.
	\end{equation}
\end{application}

\begin{application}
	[Quaternion projection on affine equality constraint] 
	The quaternion projection problem can be described as
	\begin{equation}\label{apply1}
		\begin{array}{l}
			\min\limits_{{\bm{x}}\in\mathbb{H}^n} \|\bm{x}-\bm{y}\|_{2}^{2} \\
			\text { s.t. } \bm{A x}=\bm{b}
		\end{array}
	\end{equation}
	where $\bm{y}\in \mathbb{H}^n$, $\bm{b}\in \mathbb{H}^p$, $ \bm{A} \in \mathbb{H}^{p \times n} $ and rank$(\bm{A})=p<n $. 
	Applying Example $\ref{aq-b}$, $f(\bm{x}) = \|\bm{x}-\bm{y}\|_{2}^{2}$ is convex, and $\bm{A x}=\bm{b}$ is an affine equality constraint. Therefore, the quaternion optimization problem in $(\ref{apply1})$ is convex.
	
	Using the methed of Lagrange multipliers \cite{flamant2021general,Mandic2015The}, we have
	\begin{align}
		\begin{aligned}
			&\mathcal{L}(\bm{x},\, \boldsymbol{\lambda})\\
			=
			&\|\bm{x}-\bm{y}\|_{2}^{2} + {\rm Re} \left\{\boldsymbol{\lambda}^{\mathsf{H}}(\bm{A} \bm{x}-\bm{b})\right\}\\
			=&
			(\bm{x}-\bm{y})^{\mathsf{H}} (\bm{x}-\bm{y}) + \dfrac{1}{2} \boldsymbol{\lambda}^{\mathsf{H}}(\bm{A} \bm{x}-\bm{b}) + \dfrac{1}{2} (\bm{A} \bm{x}-\bm{b})^{\mathsf{H}} \boldsymbol{\lambda}\\
			=&\bm{x}^{\mathsf{H}} \bm{x} + \left(\dfrac{1}{2} \boldsymbol{\lambda}^{\mathsf{H}}\bm{A} - \bm{y}^{\mathsf{H}}\right) \bm{x} + \bm{x}^{\mathsf{H}} \left(\dfrac{1}{2} \bm{A}^{\mathsf{H}} \boldsymbol{\lambda} - \bm{y}\right) \\
			+& \bm{y}^{\mathsf{H}} \bm{y} - \dfrac{1}{2}\boldsymbol{\lambda}^{\mathsf{H}} \bm{b} - \dfrac{1}{2} \bm{b}^{\mathsf{H}} \boldsymbol{\lambda},
		\end{aligned}
	\end{align}
	where $ \boldsymbol{\lambda} \in \mathbb{H}^{p} $ denotes the set of Lagrange multipliers. Finding the gradient of $\mathcal{L}\left(\bm{x}, \boldsymbol{\lambda}\right)$ with respect to $\boldsymbol{x}^*$ in the same way as in $(\ref{2})$ and setting the result to $\bm{0}$, we have
	\begin{align}\label{94}
		\nabla_{\bm{x}^*} \mathcal{L}(\bm{x}, \boldsymbol{\lambda})
		=&
		\left(\dfrac{\partial \mathcal{L}}{\partial \bm{x}}\right)^{\mathsf{H}}\nonumber\\
%		=&\left(\dfrac{\partial \mathcal{L}}{\partial \bm{x}}\right)^{\mathsf{H}}\nonumber\\
		=&\dfrac{1}{2} \bm{x} + \left(\dfrac{1}{2} \boldsymbol{\lambda}^{\mathsf{H}}\bm{A} - \bm{y}^{\mathsf{H}}\right)^{\mathsf{H}} - \dfrac{1}{2} \left(\dfrac{1}{2} \bm{A}^{\mathsf{H}} \boldsymbol{\lambda} - \bm{y}\right)\nonumber\\
		=&\dfrac{1}{2} \bm{x} - \dfrac{1}{2} \bm{y} + \dfrac{1}{4} \bm{A}^{\mathsf{H}} \boldsymbol{\lambda}\\
		=&\bm{0},\nonumber
	\end{align}
	which leads to
%	Setting $(\ref{94})$ to zero, and combining with $ \bm{A x}=\bm{b} $,
	\begin{align}\label{92}
		\begin{aligned}
%			\nabla_{\bm{x}^*} \mathcal{L}(\bm{x}, \boldsymbol{\lambda})
%			= \bm{0} 
%			\quad \Rightarrow \quad \ 
\bm{x}=\boldsymbol{y}-\dfrac{1}{2}  \bm{A}^{\mathsf{H}} \boldsymbol{\lambda}.
		\end{aligned}
	\end{align}
	A combinition of $(\ref{92})$ with the constraint $ \bm{A x}=\bm{b} $ yields
	\begin{align}\label{jb}
		\begin{aligned}
			&\quad\bm{A} \left(\boldsymbol{y}-\dfrac{1}{2} \bm{A}^{\mathsf{H}} \boldsymbol{\lambda}\right) = \bm{b}\\ \Rightarrow \quad &\boldsymbol{\lambda} = 2 \left(\bm{A}\bm{A}^{\mathsf{H}}\right)^{-1} (\bm{A}\bm{y}-\bm{b}). 
		\end{aligned}
	\end{align}
	Substituting $(\ref{jb})$ into $(\ref{92})$
	and applying Theorem $\ref{qjzy}$, we obtain the following optimal solution
%the closed-form optimal solution is obtained as
	\begin{equation}
		\bar{\bm{x}} = \boldsymbol{y} + \bm{A}^{\mathsf{H}} \left(\bm{A}\bm{A}^{\mathsf{H}}\right)^{-1} (\bm{b} - \bm{A}\bm{y}).
	\end{equation}
\end{application}
\begin{application}
	[Quaternion minimum variance beamforming] \label{yy1}
	%Receive Beamforming: Average Sidelobe Energy Minimization)
	%(??????: ?????????)
	The problem of quaternion variance beamforming minimization can be described as
	\begin{align}\label{y12}
		\begin{aligned}
			\min _{\bm{w} \in \mathbb{H}^n} &f(\bm{w}) \triangleq \bm{w}^{\mathsf{H}} \bm{R w} \\
			\text { s.t. } &\bm{w}^{\mathsf{H}} \bm{a}=1,
		\end{aligned}
	\end{align}
	where $\bm{w} \in \mathbb{H}^n$ is the beamformer weight vector, $\bm{a}\in \mathbb{H}^n$ is the steering vector, and $\bm{R}^{\mathsf{H}}=\bm{R}\in \mathbb{H}^{n\times n}$ is positive definite. 
%	(The first method: first-order characterization, applying Theorem $\ref{yijie}$)
	%(??2??????) 
	
	We will next prove that the problem in $(\ref{y12})$ is a convex quaternion optimization problem. 
	Using the fourth row of TABLE $\mathrm{\ref{table}}$, we have
	\begin{align}\label{111}
		\begin{aligned}
			\nabla_{\bm{w}^*} f\left(\bm{w}\right) 
			%\triangleq & \left(\dfrac{\partial f}{\partial \bm{w}^*}\right)^{\mathsf{T}} \overset{\eqref{cong}}{=} 
			=
			\left(\dfrac{\partial f}{\partial \bm{w}}\right)^{\mathsf{H}}
			= 
			\left(\boldsymbol{w}^{\mathsf{H}} \boldsymbol{R}-\frac{1}{2}\left(\boldsymbol{R} \boldsymbol{w}\right)^{\mathsf{H}}\right)^{\mathsf{H}} 
			=
			\dfrac{1}{2}\bm{R}\bm{w}.
		\end{aligned}
	\end{align}
	Then, $\forall \bm{v}, \bm{w} \in \mathbb{H}^n$,
	\begin{align}
		\begin{aligned}
				&{\rm Re}\left\{\big(\nabla_{\bm{v}^*}f\left(\bm{v}\right)-\nabla_{\bm{w}^*}f\left(\bm{w}\right)\big)^{\mathsf{H}}\left(\bm{v}-\bm{w}\right)\right\}\\
				=&\dfrac{1}{2} {\rm Re}\left\{\left(\bm{R} \bm{v} - \bm{R} \bm{w}\right)^{\mathsf{H}}\left(\bm{v}-\bm{w}\right)\right\}
				\\=&
				\dfrac{1}{2} \left(\bm{v}-\bm{w}\right)^{\mathsf{H}} \bm{R} \left(\bm{v} - \bm{w}\right)
				\\\geqslant&0.
			\end{aligned}
	\end{align}
	From Theorem $\ref{dandiao}$, 
	it follows that $f\left(\bm{w}\right)$ is convex, 
	%Applying Example $\ref{qaq}$, the quaternion function $f(\bm{w})$ is convex, while 
	and $\bm{w}^{\mathsf{H}} \bm{a}=1$ is an affine equality constraint. Therefore, the problem in $(\ref{y12})$ is a convex quaternion optimization problem.
	
	The Lagrangian of problem in $(\ref{y12})$ is given by \cite{flamant2021general,Mandic2015The}				
	\begin{equation}
		\begin{aligned}
			\mathcal{L}\left(\bm{w}, \lambda\right)  =\bm{w}^{\mathsf{H}} \bm{R} \bm{w}+\lambda\left(\bm{w}^{\mathsf{H}} \bm{a}-1\right) 
			, \quad \lambda \in \mathbb{R},
		\end{aligned}
	\end{equation}
	which is a real-valued function of $ \bm{w} \in \mathbb{H}^{n} $. %Assuming $\bm{A}$ is invertible, then we have
	Using the second and the fourth rows of TABLE $\mathrm{\ref{table}}$, and setting $\nabla_{\bm{w}^{*}} \mathcal{L}\left(\bm{w}, \lambda\right) 
	=\bm{0} $, we have
	\begin{equation}\label{115}
		\begin{aligned}
			\nabla_{\bm{w}^{*}} \mathcal{L}\left(\bm{w}, \lambda\right) 
			=
			&\left(\dfrac{\partial \mathcal{L}}{\partial \bm{w}}\right)^{\mathsf{H}}
			= \frac{1}{2}\bm{R} \bm{w}-\frac{1}{2}\lambda\bm{a}
			=\bm{0} \\
			\quad &\Rightarrow \quad \bm{w}=\lambda \bm{R}^{-1} \bm{a} .
		\end{aligned}
	\end{equation}				
	Upon substituting $(\ref{115})$ into $ \bm{a}^{\mathsf{H}} \bm{w}=1 $, we obtain
	
	\begin{equation}
		\lambda \bm{a}^{\mathsf{H}} \bm{R}^{-1} \bm{a}=1 \quad  \Rightarrow \quad  \lambda=\frac{1}{\bm{a}^{\mathsf{H}} \bm{R}^{-1} \bm{a}} .
	\end{equation}
	Therefore, upon applying Theorem $\ref{qjzy}$, the closed-form optimal solution is obtained as
	\begin{equation}
		\bar{\bm{w}}=\frac{\bm{R}^{-1} \bm{a}}{\bm{a}^{\mathsf{H}} \bm{R}^{-1} \bm{a}} .
	\end{equation}

\end{application}
\section{Conclusions}
%\textcolor{blue}{
	We have established the theory of convex quaternion optimization based on the GHR calculus, which is an enabling methodology in the field of quaternion optimization and its applications in quaternion signal processing and machine learning. Our study has resulted in the development of five discriminant theorems for convex functions in the quaternion field, utilizing $(\ref{r,h})$, $(\ref{h})$, $(\ref{HHH})$, $(\ref{grad})$, and  $(\ref{hassion})$. Furthermore, we have provided the definition and four discriminant criteria for strongly convex functions by employing the results for convex quaternion functions.
	In addition, we have presented a fundamental theorem for the optimality of convex quaternion optimization problems and three applications in signal processing, which have both enriched the theory of convex quaternion optimization and provided a theoretical foundation for quaternion signal processing. However, the convexity of non-differentiable quaternion functions by the GHR calculus still remains an open area, and this work provides a foundation and an avenue for further research in this direction.

\bibliographystyle{IEEEtran}
\bibliography{IEEEabrv,bib6}

\end{document}